\documentclass[11pt, twoside]{amsart}

\usepackage[utf8]{inputenc}
\usepackage[T1]{fontenc}
\usepackage{amssymb,amsmath,amstext, graphicx, hyperref}

\newtheorem{theor}{Theorem}[section]
\newtheorem{lem}[theor]{Lemma}
\newtheorem{defin}[theor]{Definition}

\newtheorem{prop}[theor]{Proposition} 

\newtheorem{notation}[theor]{Notation}
\newtheorem{exam}[theor]{Example}

\newtheorem{cor}[theor]{Corollary}
\newtheorem{rem}[theor]{Remark}

\newtheorem{fact}[theor]{Fact}

\newtheorem{assump}[theor]{Assumption}
\newtheorem{observation}[theor]{Observation}

\numberwithin{equation}{section}

\newcommand{\tp}{\mathrm{tp}}

\newcommand{\acl}{\mathrm{acl}}
\newcommand{\dcl}{\mathrm{dcl}}

\newcommand{\cl}{\mathrm{cl}}

\newcommand{\mb}{\mathbf}

\newcommand{\es}{\emptyset}

\newcommand{\su}{\mathrm{SU}}

\newcommand{\nts}{\negthickspace}
\newcommand{\uhrc}{\nts \upharpoonright \nts}

\newcommand{\crd}{\mathrm{crd}}

\newcommand{\meq}{^{\mathrm{eq}}}

\newcommand{\mcA}{\mathcal{A}}
\newcommand{\mcB}{\mathcal{B}}
\newcommand{\mcC}{\mathcal{C}}
\newcommand{\mcD}{\mathcal{D}}

\newcommand{\mcG}{\mathcal{G}}

\newcommand{\mcM}{\mathcal{M}}
\newcommand{\mcN}{\mathcal{N}}

\newcommand{\mcR}{\mathcal{R}}

\newcommand{\mbP}{\mathbf{P}}

\newcommand{\mbK}{\mathbf{K}}

\newcommand{\real}{\mathrm{real}}

\newcommand{\ind}{\raisebox{-2pt}[5pt][0pt]{$\smile$} \hspace*{-6.8pt}\raisebox{3pt}[5pt][0pt]{$|$} \; \: }
\newcommand{\nind}{\raisebox{-2pt}[5pt][0pt]{$\smile$} 
\hspace*{-6.8pt}\raisebox{3pt}[5pt][0pt]{$|$}\hspace*{-6.8pt}
\raisebox{3pt}[5pt][0pt]{$\diagup$} }

\newcommand{\rng}{\mathrm{rng}}

\title[Binary primitive homogeneous simple structures]
{Binary primitive homogeneous simple structures}

\author{Vera Koponen}

\address{Vera Koponen, Department of Mathematics, Uppsala University, Box 480,
75106 Uppsala, Sweden.}

\email{vera.koponen@math.uu.se}

\date{22 July 2016 (revised version)}

\begin{document}

\begin{abstract}
Suppose that $\mcM$ is countable, binary, primitive, homogeneous, and simple, and hence 1-based. We prove that the SU-rank of the complete theory of $\mcM$ is~1. It follows that $\mcM$ is a random structure. The conclusion that $\mcM$ is a random structure does not hold if the binarity condition is removed, as witnessed by the generic tetrahedron-free 3-hypergraph. However, to show that the generic tetrahedron-free 3-hypergraph is 1-based requires some work (it is known that it has the other properties) since this notion is defined in terms of imaginary elements. This is partly why we also characterize equivalence relations which are definable without parameters in the context of $\omega$-categorical structures with degenerate algebraic closure.  Another reason is that such characterizations may be useful in future research about simple (nonbinary) homogeneous structures.

\end{abstract}

\maketitle

\section{Introduction}\label{Introduction}

\noindent
This article is part of a study of (in particular binary) homogeneous and simple structures.
In order not to be too repetitive we refer to the introductory sections of \cite{AK, Kop14} 
for more background concerning homogeneous structures,
simple structures and the conjunction of both.
But in general the idea is that, although some particular classes of homogeneous structures have been classified,
the class of all homogeneous structures is too large and diverse to be understood in a relatively uniform way.\footnote{
For a survey of homogeneous structures, 
including applications to permutation groups, constraint satisfaction problems, Ramsey theory and topological dynamics,
see~\cite{Mac10}.
For a classification of all homogeneous digraphs, see~\cite{Che98}. Both sources contain
many references, for example to classifications of various kinds of homogeneous structures.}
So we like to impose some extra conditions that give us tools to work with.
Given the existing and model theoretically important framework of simple structures 
\cite{Wag} it is natural to consider structures which are both homogeneous and simple.
The probably most well known example in this class is the Rado graph, an example of a random structure in the
sense of Definition~\ref{definition of a random structure} below.
The study of homogeneous simple structures is also an extension of the work of Lachlan and others about
stable homogeneous structures~\cite{Lach97}.

Here a structure $\mcM$ is called {\em homogeneous} if it has a finite relational vocabulary (signature) and 
every isomorphism between finite substructures can be extended to an automorphism of $\mcM$.
If $\mcM$ is a countable structure with finite relational vocabulary, then $\mcM$ is homogeneous if and only
if it has elimination of quantifiers. 
When assuming that a structure is simple we automatically assume that it is infinite.
A structure $\mcM$ is called {\em primitive} if there is no {\em nontrivial} equivalence relation on its universe $M$
which is $\es$-definable, i.e. definable without parameters. (By a nontrivial equivalence relation we mean one
which has at least two equivalence classes and at least one equivalence class contains more than one element.)
A reason why primitive homogeneous structures are of interest is the following: Suppose that $\mcM$ is a
homogeneous structure with a nontrivial $\es$-definable equivalence relation on $M$. 
Let $A$ be any one of the equivalence classes. Then it is easy to see that $\mcM \uhrc A$, the
substructure of $\mcM$ with universe $A$, is homogeneous. If $\mcM$ is, in addition, simple and $A$ is infinite,
then $\mcM \uhrc A$ is also simple. Thus we cannot understand $\mcM$ any better than we can understand 
$\mcM \uhrc A$. If, in particular, $A$ is an equivalence class which cannot be ``split'' into two nonempty parts by
some other $\es$-definable equivalence relation, then $\mcM \uhrc A$ is a primitive structure.

A structure $\mcM$ is called {\em binary} if its vocabulary contains only unary and/or binary relation symbols.
For basics about simple structures see for example~\cite{Wag}.
Our first (and main) result is the following:

\begin{theor}\label{main theorem}
Suppose that $\mcM$ is countable, binary, primitive, homogeneous and simple.
Then the SU-rank of $Th(\mcM)$ is 1.
\end{theor}

\noindent
By a result of Aranda Lopez \cite[Proposition~3.3.3]{AL}, if
$\mcM$ is binary, homogeneous, primitive and simple and $Th(\mcM)$ has SU-rank 1,
then $\mcM$ is a random structure in the sense of 
Definition~\ref{definition of a random structure} below.\footnote{
Proposition 3.3.3 in \cite{AL} does not use the terminology ``random structure'' but formulates the result
in terms of ``Alice's restaurant property'',
also known as ``extension axioms/properties''.}
Hence we get the following consequence, which gives a positive answer to the leading question asked by the author in \cite{Kop14}:

\begin{cor}\label{corollary to main theorem}
If $\mcM$ is countable, binary, primitive, homogeneous and simple, then $\mcM$ is a random structure.
\end{cor}

\noindent
As shown by Example~\ref{binarity cannot be removed},
the binarity assumption cannot be removed from Theorem~\ref{main theorem}.
Since random structures (according to Definition~\ref{definition of a random structure}) have SU-rank 1, it follows that
the binarity condition cannot be removed from Corollary~\ref{corollary to main theorem}.

The proof of Theorem~\ref{main theorem} relies on the following claim:
{\em If $\mcM$ is countable, binary, homogeneous and simple, then it is 1-based.}
This is also stated as (part of) Fact~\ref{binary homogeneous simple structures are 1-based} below, and a justification
for this claim is given in Remark~\ref{proof of an earlier fact}.
Actually, all currently known homogeneous simple structures are 1-based, or appear to be so.
(We will see that verifying 1-basedness is not necessarily straightforward, even if the dividing/forking behaviour on real elements
is as simple as it can be.)

Since the binarity condition cannot be removed from Corollary~\ref{corollary to main theorem} one may ask the following question:
If $\mcM$ is countable, primitive, homogeneous, supersimple with SU-rank 1 and 1-based, must $\mcM$ be a random structure?
The answer is no, as witnessed by the generic tetrahedron-free 3-hypergraph 
from Definition~\ref{definition of generic tetrahedron-free 3-hypergraph} below:

\begin{prop}\label{the generic tetrahedron free 3-hypergraph is a counterexample}
The generic tetrahedron-free 3-hypergraph 
is primitive, homogeneous, supersimple with SU-rank 1 and 1-based, but not a random structure.
\end{prop}

\noindent
It is known that the generic tetrahedron-free 3-hypergraph is primitive, homogeneous (by construction)
and supersimple with SU-rank 1, but not a random structure 
(in the sense of Definition~\ref{definition of a random structure} below). 
See 
Remark~\ref{remark about the generic tetrahedron-free 3-hypergraph}
for further explanations. However, as far as the author knows, the claim of
Proposition~\ref{the generic tetrahedron free 3-hypergraph is a counterexample}
that it is 1-based has never been verified before
and is not a trivial matter, because we need to deal with imaginary elements, in other words with elements defined by
$\es$-definable equivalence relations on tuples of elements from the structure.
(Independently of the present author, Conant has recently proved a result, about homogeneous structures whose age
has ``free amalgamation'', which implies that the 
generic tetrahedron-free 3-hypergraph is 1-based \cite{Con}.)

Thus we prove two results, Theorems~\ref{characterization of equivalence relations for omega-categorical} 
and~\ref{characterization of equivalence relations in a special case}, 
which characterize $\es$-definable equivalence relations on 
tuples of elements, starting from a bit different assumptions.
Besides being used to justify the 1-basedness claim of 
Proposition~\ref{the generic tetrahedron free 3-hypergraph is a counterexample}, via
Example~\ref{examples satisfying the special case}~(ii) and
Proposition~\ref{the generic tetrahedron-free hypergraph is 1-based} below,
these theorems may be useful in the future for understanding $\omega$-categorical structures with additional 
properties, such as being homogeneous and simple (but not necessarily binary).
Corollary~\ref{algebraic closure equals definable closure} shows that,
under the same hypotheses as in
Theorem~\ref{characterization of equivalence relations in a special case},
the algebraic closure and definable closure in $\mcM\meq$ of any $A \subseteq M$ are identical.

Now follows an outline of this article and of the proof of Theorem~\ref{main theorem}.
Section~\ref{Preliminaries} 
gives a few definitions and remarks of relevance for this article. Since the work here takes place
within the same context as \cite{AK, Kop14} we refer to the preliminary section of any one of these articles
for more detailed explanations of notions and known results 
that will be used (concerning homogeneous and $\omega$-categorical structures with simple theories and about imaginary elements).

Suppose that $\mcM$ is countable, binary, primitive, homogeneous and simple.
Then $\mcM$ is supersimple with finite SU-rank, 1-based and has trivial dependence.
(See Fact~\ref{binary homogeneous simple structures are 1-based} and the discussion just before and after it.)
Hence the results about coordinatization developed in \cite[Section~3]{Djo06} are applicable to $\mcM$.
These results and \cite[Theorem~5.1]{AK} were used in \cite{Kop14} to show that
$\mcM$ can be ``strongly interpreted'' in a binary random structure. This ``strong interpretation'' can also
be seen as a coordinatization of  $\mcM$ by a binary random structure and constitutes the framework within which
we will prove Theorem~\ref{main theorem}. This framework is explained in 
Section~\ref{coordinatization of M}.
In Section~\ref{Proof of the main theorem} we prove Theorem~\ref{main theorem}.

As already mentioned, the main results of Sections~\ref{Definable equivalence relations}
and~\ref{A variation on the theme}  
(Theorems~\ref{characterization of equivalence relations for omega-categorical} 
and~\ref{characterization of equivalence relations in a special case})
characterize $\es$-definable equivalence relations on $n$-tuples
($0 < n < \omega$) under assumptions including $\omega$-categoricity and degenerate algebraic closure.
These sections do not depend on
Sections~\ref{coordinatization of M} or~\ref{Proof of the main theorem} 
and can be read separately.
Theorem~\ref{characterization of equivalence relations in a special case}
is used to prove Proposition~\ref{the generic tetrahedron free 3-hypergraph is a counterexample},
via Proposition~\ref{the generic tetrahedron-free hypergraph is 1-based}.

\section{Preliminaries}\label{Preliminaries}

\noindent
The notation and terminology used here is more or less standard, but we nevertheless begin with
clarifying some notation. First-order structures (the only kind considered) are denoted
$\mcA, \mcB, \ldots, \mcM, \mcN, \ldots$ and their universes are denoted $A, B, \ldots, M, N, \ldots$, respectively.
Finite sequences are denoted by $\bar{a}, \bar{b}, \ldots, \bar{x}, \bar{y}, \ldots$.
We may denote the concatenation of $\bar{a}$ and $\bar{b}$ by $\bar{a}\bar{b}$.
The set of elements occuring in $\bar{a}$ is denoted by $\rng(\bar{a})$, ``the range of $\bar{a}$''.
We often write `$\bar{a} \in A$' as shorthand for `$\rng(\bar{a}) \subseteq A$'.
A structure is called $\omega$-categorical, (super)simple or 1-based if its complete theory has the corresponding property.
The SU-rank of a supersimple structure is (by definition) the SU-rank of its complete theory;
and the SU-rank of a supersimple complete theory $T$ is the supremum (if it exists) of the SU-ranks of all
1-types over $\es$ with respect to $T$. 
In this article it is often important to distinguish in which structure a complete type, the algebraic closure
etcetera, is taken, so we use subscripts (or superscripts) such as in `$\tp_\mcM$' or `$\acl_{\mcM\meq}$' to indicate this.
It will also be convenient to occasionally use the notation
$\bar{a} \equiv_\mcM \bar{b}$ as shorthand for $\tp_\mcM(\bar{a}) = \tp_\mcM(\bar{b})$.

The context of this article is the same as that of \cite{AK, Kop14}  and therefore we refer to those
articles (any one of them will do) for basics and relevant facts about
homogeneous structures, simple structures and the extension $\mcM\meq$ of $\mcM$ by imaginaries.
However we repeat the following definitions here:
we say that $\mcN$ is {\em canonically embedded} in $\mcM\meq$
if $N$ is a $\es$-definable subset of $M\meq$ and for every $0 < n < \omega$ and every 
relation $R \subseteq N^n$ that is $\es$-definable in $\mcM\meq$ there is a relation symbol in the
vocabulary of $\mcN$ which is interpreted as $R$, and the vocabulary of $R$ contains no other symbols.
If $\mcM$ and $\mcN$ are structures (possibly with different vocabularies) then $\mcN$ is a {\em reduct}
of $\mcM$ if the following holds: $M = N$ and if $0 < k < \omega$ and $R \subseteq N^k$ is $\es$-definable in $\mcN$,
then $R$ is $\es$-definable in $\mcM$.

Also, it is important to distinguish between two distinct, but related, notions of ``triviality''.
A pregeometry (or matroid, see \cite[Chapter~4.6]{Hod} for a definition) $(A, \cl)$ will be called {\em trivial} if, 
for all $a \in A$ and $B \subseteq A$,
$a \in \cl(B)$ implies that $a \in \cl(\{b\})$ for some $b \in B$.
A structure $\mcM$ has {\em degenerate algebraic closure} if for every $A \subseteq M$, $\acl_\mcM(A) = A$;
in this case we may also say that $\acl_\mcM$ is {\em degenerate}.

\noindent
The question of what a ``truly''  random structure is does not have an obvious answer, but
here is the definition that we will use:

\begin{defin}\label{definition of a random structure}{\rm
(i) Let $V$ be a vocabulary and let $\mcM$ be a $V$-structure.
We call a finite $V$-structure $\mcA$ a {\em forbidden structure with respect to $\mcM$} if $\mcA$ cannot be
embedded into $\mcM$. If, in addition, there is no proper substructure of $\mcA$ which is forbidden with respect to $\mcM$,
then we call $\mcA$ a {\em minimal forbidden structure with respect to $\mcM$}.\\
(ii) If $W \subseteq V$ are vocabularies and $\mcM$ is a $V$-structure, then $\mcM \uhrc W$ denotes the reduct
of $\mcM$ to $W$.\\
(iii) Let $V$ be a finite relational vocabulary with maximal arity $r$, where $r \geq 2$. 
We say that a $V$-structure $\mcM$ is a {\em random structure} if $\mcM$ is infinite, countable, homogeneous and, 
for every $k = 2, \ldots, r$,
there does {\em not} exist a minimal forbidden structure $\mcA$ with respect to 
\[
\mcM \uhrc \{ P \in V : \text{ the arity of $P$ is } \leq k\}
\]
such that $|A| \geq k+1$.
If $\mcM$ is a random structure and the maximal arity of 
its vocabulary is 2, then we may call $\mcM$ a {\em binary random structure}.
}\end{defin}

\begin{rem}\label{remark about random structures}{\rm
The definition of {\em binary random structure} above coincides with the one given in \cite{Kop14} and
is equivalent to the definition given in \cite{AK}.
Clearly, the Rado graph is a binary random structure according to the definition given here.
}\end{rem}

\begin{defin}\label{definition of generic tetrahedron-free 3-hypergraph}{\rm
(i) A {\em 3-hypergraph} is a structure $\mcM$ whose vocabulary contains one ternary relation symbol, say $P$,
(and no other symbols) and which satisfies the following for any permutation $\pi$ of $\{1, 2, 3\}$:
\[
\forall x_1, x_2, x_3 \Bigg( P(x_1, x_2, x_3) \ \rightarrow \ \bigg[\bigwedge_{i \neq j} x_i \neq x_j \ \wedge \
P(x_{\pi(1)}, x_{\pi(2)}, x_{\pi(3)}) \bigg] \Bigg).
\]
(ii) By a {\em tetrahedron} we mean a 3-hypergraph $\mcM$ such that $|M| = 4$ and for all distinct $a, b, c \in M$, 
$\mcM \models P(a, b, c)$. A 3-hypergraph $\mcM$ is called {\em tetrahedron-free} if no tetrahedron can be embedded into it.\\
(ii) Let $\mbK$ be the class of all finite tetrahedron-free 3-hypergraphs. Then $\mbK$ has the hereditary property and
amalgamation property and therefore $\mbK$ has a (unique) Fra\"{i}ss\'{e} limit, which is an infinite countable homogeneous
structure.\footnote{
See \cite[Chapter~7]{Hod} for the involved notions and relevant results. 
It is straightforward to see that $\mbK$ has the hereditary 
property and amalgamation property.
The joint embedding property follows from the amalgamation property since the
vocabulary is relational.}
}\end{defin}

\begin{rem}\label{remark about the generic tetrahedron-free 3-hypergraph}{\rm
Let $\mcM$ be the generic tetrahedron-free 3-hypergraph.
It is known that $\mcM$ is supersimple with SU-rank 1 and an argument showing this is found in
\cite[Section~3]{Djo06a} where the same structure is called the ``random pyramid-free (3)-hypergraph''.
Moreover (see~\cite{Djo06a}), if $A, B, C  \subseteq M$ then  $A \underset{C}{\ind} B$ if and only if
$A \cap (B \cup C) = A \cap C$.
Since all pairs of distinct elements have the same type it follows that $\mcM$ is also primitive.
Clearly, $\mcM$ is not a random structure.
Proposition~\ref{the generic tetrahedron-free hypergraph is 1-based} implies that $\mcM$ is 1-based. Its proof uses 
Theorem~\ref{characterization of equivalence relations in a special case}, the proof of which is a slight variation of the proof of 
Theorem~\ref{characterization of equivalence relations for omega-categorical}.
}\end{rem}

\noindent
The useful consequence of 1-basedness in the context of homogeneous simple structures is that dependence is trivial,
in the sense of the following definition:

\begin{defin}\label{definition of trivial dependence}{\rm
Let $\mcM$ be a simple structure. 
We say that $\mcM$ has {\em trivial dependence} if whenever $\mcN \models Th(\mcM)$,
$A, B, C_1, C_2 \subseteq N\meq$ and $A \underset{B}{\nind} (C_1 \cup C_2)$, then 
$A \underset{B}{\nind} C_i$ for $i = 1$ or $i = 2$.
}\end{defin}

If $\mcM$ is countable, homogeneous and simple, then (by results of Macpherson~\cite{Mac91}, De Piro and Kim~\cite{PK}
and Hart, Kim and Pillay~\cite{HKP})
$\mcM$ is 1-based if and only if it has trivial dependence and finite SU-rank.
By \cite[Theorem 1]{Kop15}, every countable, binary, homogeneous simple structure is supersimple with finite SU-rank.
Therefore, if $\mcM$ is countable, binary, homogeneous and simple, then it is 1-based if and only if it has trivial dependence.

\begin{fact}\label{binary homogeneous simple structures are 1-based}
If $\mcM$ is countable, binary, homogeneous and simple then it is supersimple with finite SU-rank, 1-based and has trivial dependence.
\end{fact}

\noindent
From what has been said it follows that the only thing that needs to be proved is that $\mcM$, 
as in Fact~\ref{binary homogeneous simple structures are 1-based}, has trivial dependence.
We postpone this to Remark~\ref{proof of an earlier fact} at the very end,
because then we can ``reuse'' a part of the argument in the proof of 
Proposition~\ref{the generic tetrahedron-free hypergraph is 1-based}, rather than repeating that argument.

\begin{exam}\label{binarity cannot be removed}{\rm
The necessity of the binarity assumption in Theorem~\ref{main theorem} is shown by the following example which 
also appears as Example~3.3.2 in \cite{Mac10}.

Let $\mcM$ be a countable infinite structure with empty vocabulary, so $\mcM$ is just a set, 
hence $\mcM$ is homogeneous and $\omega$-stable (thus supersimple) of SU-rank 1. 
Trivially, $\mcM$ is binary, so by Fact~\ref{binary homogeneous simple structures are 1-based} $\mcM$ is 1-based
and has trivial dependence.
Let $\mcG = (V, E)$ be the graph where $V$ is the set of all (unordered) 2-subsets of $M$ and let two vertices of $G$ 
be adjacent if and only if they intersect in exactly one point (of $M$). Then $\mcG$ is interpretable in $M$ (without parameters), 
so it is $\omega$-categorical and stable. 
However $\mcG$ is not homogeneous. Let $\mcG'$ be the expansion of $\mcG$ by adding a ternary relation symbol $Q$, 
where  $Q(a, b, c) $ holds in  $\mcG'$ if and only if  $a$, $b$ and $c$ are distinct and the intersection of all three is nonempty. 
Then $\mcG'$ is homogeneous (we leave the proof to the reader). 
Moreover, $Q$ is definable in $\mcG$ without parameters, so $\mcG'$ is stable. 
Observe that every permutation of $M$ naturally induces an automorphism of $\mcG'$
(if distinct $a, b \in M$ are mapped to $a', b'$, respectively, then let $\{a, b\}$ be mapped to $\{a', b'\}$).
Therefore $\mcG'$ has a unique 1-type over $\es$ and there are exactly 2 different 2-types of distinct elements over $\es$
(adjacent or nonadjacent vertices).
By using the definition of dividing it is straightforward to show that the unique 1-type (over $\es$) of $\mcG'$ has SU-rank 2, 
so $\mcG'$ has SU-rank 2.\footnote{
The following argument shows that the SU-rank is at least 2.
Let $a, b \in V$ and $E(a, b)$. Then $\tp_{\mcG'}(a /b)$ is nonalgebraic so $\su(a / b) \geq 1$.
Take $b_i \in V$, $i < \omega$, such that $b = b_0$ and whenever $i \neq j$ then $\neg E(b_i, b_j)$.
Then $\{E(x, b_i) : i < \omega\}$ is 3-inconsistent, so $\tp_{\mcG'}(a / b)$ divides over $\es$. Hence $\su(a) \geq 2$.}

The vocabulary of $\mcG'$ has only the symbols $E$ and $Q$ and it is easy to see that $E$ is not an equivalence relation.
As $\mcG'$ has elimination of quantifiers (being homogeneous) it follows that it is primitive.
Furthermore, $\mcG'$ is 1-based. To show this, it is, by \cite[Corollary 4.7]{HKP}, sufficient to show that
{\em for every complete type of SU-rank 1 (possibly realized by imaginary elements),
the pregeometry on its realizations (given by algebraic closure) is trivial.}
Since $\mcM$ has trivial dependence this is true for $\mcM$. 
As $\mcG'$ is definable without parameters in $\mcM\meq$ it follows that 
$\mcG'^{\mathrm{eq}}$ is definable without parameters in $(\mcM\meq)\meq$. 
Since $\mcM\meq$ has elimination of imaginaries it follows that $\mcG'^{\mathrm{eq}}$ is definable without parameters
in $\mcM\meq$.  Hence the statement in italics holds for $\mcG'$.\footnote{
Strictly speaking, the statement in italics should be proved in the context when the parameters of the type
come from an arbitrary model of $Th(\mcG')$, so the verification of 1-basedness is not quite complete.
However, the beginning of the proof of
Proposition~\ref{the generic tetrahedron-free hypergraph is 1-based} 
shows how to overcome this slight obstacle.}
}\end{exam}

\section{Coordinatization by a random structure}
\label{coordinatization of M}

\noindent
In this section and in Section~\ref{Proof of the main theorem} 
we assume that $\mcM$ {\bf \em is countable, binary, primitive, homogeneous and simple}.
By Fact~\ref{binary homogeneous simple structures are 1-based},
{\bf \em $\mcM$ is supersimple with finite SU-rank, 1-based and has trivial dependence.
Let the SU-rank of $\mcM$ be $\rho$.}

Note that the primitivity of $\mcM$ implies that for all $a, b \in M$, $\tp_\mcM(a) = \tp_\mcM(b)$.

\begin{fact}\label{existence of C} (\cite[Section~3]{Kop14}, which uses the ``coordinatization'' from \cite{Djo06})
There is $C \subseteq M\meq$ such that:
\begin{itemize}
\item[(i)] $C$ is $\es$-definable in $\mcM\meq$ and only finitely many sorts are represented in $C$,
\item[(ii)] $\su(c) = 1$ for every $c \in C$ (where SU-rank of elements/types is taken with respect to $\mcM\meq$),
\item[(iii)] $\acl_{\mcM\meq}$ restricted to $C$ is {\em degenerate}, by which we mean that
$\acl_{\mcM\meq}(A) \cap C = A$ for every $A \subseteq C$, and
\item[(iv)] for every $a \in M$, $a \in \acl_{\mcM\meq}(\crd(a))$ where we define
$\crd(a) = \acl_{\mcM\meq}(a) \cap C$
(so in particular $M \subseteq \acl_{\mcM\meq}(C)$).
\item[(v)] For every $c \in C$ there is $a \in M$ such that $c \in \crd(a)$.\footnote{
This is the only part which may not be immediate from \cite[Section~3]{Kop14}.
However, if $C$ has all properties~(i)--(iv) but not~(v), then
we can let $C' = \{ c \in C : c \in \crd(a) \text{ for some } a \in M\}$ and it is straightforward to verify
(using that $\mcM$ is $\omega$-categorical)
that $C'$ satisfies~(i)--(v), because for every $a \in M$, $\crd(a)$ is the same whether computed with respect 
to $C$ or with respect to $C'$.}
\end{itemize}
\end{fact}

\noindent
We call $C$ as in Fact~\ref{existence of C}
a set of {\em coordinates} of $\mcM$ and for each $a \in M$, $\crd(a)$ may
be called the (set of) {\em coordinates of $a$}.
Moreover, if $A \subseteq C$ and there is $a \in M$ such that $\crd(a) = A$, then we call $A$ a {\em line}.
From the assumptions about $\mcM$ and the properties of $C$ and $\crd$ from
Fact~\ref{existence of C}
one easily derives the following\footnote{For part (ii), use \cite[Lemma~3.5]{Kop14}, primitivity of $\mcM$,
and the fact that (by $\omega$-categoricity) the equivalence relation `$\acl_{\mcM}(x) = \acl_{\mcM}(y)$' 
has infinitely many classes}:

\begin{fact}\label{derived properties of C and crd}
Let $C$ and $\crd$ be as in Fact~\ref{existence of C}.
Then:
\begin{itemize}
\item[(i)] For all $a, a' \in M$, $\su(a) = |\crd(a)| = |\crd(a')| = \su(a')$, so in particular all lines have the
same cardinality, which is $\rho$.
\item[(ii)] For all $a \in M$, $a \in \dcl_{\mcM\meq}(\crd(a))$.
Hence, for every line $A \subseteq C$ there is a unique $a \in M$ such that $\crd(a) = A$.
\item[(iii)] for all $a, a' \in M$, $\crd(a)$ and $\crd(a')$ can be ordered as $\bar{c}$ and $\bar{c}'$, respectively,
so that $\tp_{\mcM\meq}(\bar{c}) = \tp_{\mcM\meq}(\bar{c}')$.
\end{itemize}
\end{fact}

\begin{assump}\label{assumption about C}
For the rest of this section and in Section~\ref{Proof of the main theorem} we assume that 
{\bf \em $C$ is a coordinatizing set as in Fact~\ref{existence of C}
and $\mcC$ is the canonically embedded structure in $\mcM\meq$ with universe $C$.}
\end{assump}

\noindent
From \cite[Theorem~5.1]{AK} we immediately get the following:

\begin{fact}\label{C is a reduct of a binary random structure}
There is a binary random structure $\mcR$ such that $\mcC$ is a reduct of $\mcR$.
\end{fact}

\noindent
Moreover, by \cite[Lemma~3.9]{Kop14}:

\begin{fact}\label{lines are rigid in R}
Suppose that $\mcR$ is like in Fact~\ref{C is a reduct of a binary random structure}.
Let $a \in M$ and $\crd(a) = \{c_1, \ldots, c_\rho\}$ where the elements are enumerated without repetition.
Then for every nontrivial permutation $\pi$ of $\{1, \ldots, \rho\}$,
$(c_1, \ldots, c_\rho) \not\equiv_\mcR (c_{\pi(1)}, \ldots, c_{\pi(\rho)})$.
\end{fact}

\noindent
Unfortunately, the information given by these facts is not quite enough for the purpose of proving 
Theorem~\ref{main theorem}.
Therefore the next result gives the strengthening of Fact~\ref{C is a reduct of a binary random structure} that we need.

\begin{fact}\label{getting suitable R}
There is a binary random structure $\mcR$ such that:
\begin{itemize}
\item[(i)] $\mcC$ is a reduct of $\mcR$.

\item[(ii)] For every $p(x) \in S^{\mcM\meq}_1(\acl_{\mcM\meq}(\es))$
that is realized in $C$ there is a unique unary relation symbol $R_p$ in the vocabulary of $\mcR$ such that
for all $c \in C$, $\mcR \models R_p(c)$ if and only if $\mcM\meq \models p(c)$.

\item[(iii)] For every $q(x, y) \in S^{\mcM\meq}_2(\es)$ that is realized in $C^2$
there is a unique binary relation symbol $R_q$ in the vocabulary of $\mcR$ such that for all $c_1, c_2 \in C$,
$\mcR \models R_q(c_1, c_2)$ if and only if $\mcM\meq \models q(c_1, c_2)$.

\item[(iv)] The vocabulary of $\mcR$ has no other symbols than those mentioned in~(ii) and~(iii).
Consequently, every type in $S^\mcR_n(\es)$ (for any $0 < n <  \omega$) is isolated by a
conjunction of such formulas $R_p(x) $ and $R_q(x,y)$ mentioned in~(ii) and~(iii).
\end{itemize}
\end{fact}

\noindent
{\em Remark:} In~(i) we consider 1-types over $\acl_{\mcM\meq}(\es)$ and in~(ii)
we consider 2-types over $\es$; this is the intention and not a mistake.
\\

\noindent
{\bf Proof.}
The proof is a modification of the proof of Theorem~5.1 in~\cite{AK} and we only explain
how to modify that proof.
By Lemma~4.5 in~\cite{AK} there is a $\es$-definable $D \subseteq M\meq$ 
in which only finitely many sorts are represented and
such that 
\begin{itemize}
\item[(a)] $D$ is $\acl$-complete (see Definition~4.4 in~\cite{AK}),

\item[(b)] for every $c \in C$ there is $d \in D$ such that $c \in \dcl_{\mcM\meq}(d)$ and $d \in \acl_{\mcM\meq}(c)$, and

\item[(c)] for every $d \in D$ there is (a unique) $c \in C$ such that 
$c \in \dcl_{\mcM\meq}(d)$ and $d \in \acl_{\mcM\meq}(c)$.
\end{itemize}
Note that if $c$ and $d$ are as in~(b) (or~(c)), then $\tp_{\mcM\meq}(c / \acl_{\mcM\meq}(\es))$
is determined by $\tp_{\mcM\meq}(d / \acl_{\mcM\meq}(\es))$, because 
$c \in \dcl_{\mcM\meq}(d)$. If, in addition, $c'$ and $d'$ satisfy the
same conditions as $c$ and $d$, then $\tp_{\mcM\meq}(c, c')$ is determined by $\tp_{\mcM\meq}(d, d')$.

Just as in the proof of Theorem~5.1 in~\cite[p 244]{AK}, let 
$p_1, \ldots, p_r$ be all complete 1-types over $\acl_{\mcM\meq}(\es)$ which 
are realized in $D$, and let $p_{r+1}, \ldots, p_s$  be all complete 2-types over $\es$
which are realized in $D^2$. 
For each $i = 1, \ldots, r$ let $R_i$ be a unary relation symbol and for each 
$i = r+1, \ldots, s$ let $R_i$ be a binary relation symbol.
Let $V = \{R_1, \ldots, R_s\}$ and let $\mcD$ be the $V$-structure with universe $D$ such that,
for every $i = 1, \ldots, s$ and every $\bar{d} \in D$ of appropriate length,
$\mcD \models R_i(\bar{d})$ if and only if $\mcM\meq \models p_i(\bar{d})$.

So far we have followed the proof of Theorem~5.1 in~\cite{AK}.
The difference comes now when we define a subvocabulary $V' \subseteq V$ and then a class
$\mbK$ of finite $V'$-structures, instead of a class $\mbK$ of finite $V$-structures as in~\cite{AK}.
Let $I \subseteq \{1, \ldots, s\}$ be minimal (with respect to inclusion) such that the following hold:
\begin{itemize}
\item For every $q \in S^{\mcM\meq}_1(\acl_{\mcM\meq}(\es))$ that is realized in $C$ there is 
$p_i$, $1 \leq i \leq r$ and $i \in I$,
such that whenever $c \in C$ and $\mcM\meq \models q(c)$, then there is $d \in D$ satisfying~(b)
and $\mcM\meq \models p_i(d)$.

\item For every $q \in S^{\mcM\meq}_2(\es)$ that is realized in $C^2$ there is 
$p_i$, $r < i \leq s$ and $i \in I$, such that whenever $c_1, c_2 \in C$ and $\mcM\meq \models q(c_1, c_2)$,
then there are $d_1, d_2 \in D$ satisfying~(b) ($d_i$ with respect to $c_i$) and
$\mcM\meq \models p_i(d_1, d_2)$.
\end{itemize}

\noindent
Now let
$
V' = \{R_i : i \in I\}
$
and let $\mbK$ be the class of all finite $V'$-structures $\mcN$ such that
\begin{itemize}
\item for all $a \in N$, $\mcN \models R_i(a)$ for some $R_i \in V'$, and

\item there is an embedding $f: \mcN \to \mcD \uhrc V'$ such that
$f(N)$ is an independent set (where independence is with respect to $\mcM\meq$).
\end{itemize}

\noindent
Now we can define $\mbP_2$ and $\mb{RP}_2$ from $\mbK$ in the same way as in \cite[p 244]{AK}.
The rest of the proof, starting from Lemma~5.3 in~\cite{AK}, is like the proof of Theorem~5.1 in~\cite{AK};
although $\mbK$ (and consequently $\mbP_2$ and $\mb{RP}_2$) is defined differently here the 
same arguments work out in the present context.
Hence we find a binary random $V'$-structure $\mcR$ such that~(i)--(iii) of this lemma hold.
This modification of the proof of \cite[Theorem~5.1]{AK} thus amounts to showing that the vocabulary
$V$ may have redundant symbols (for the purpose of making $\mcC$ a reduct of $\mcR$)
and we can always do with the vocabulary $V'$ as defined above.
\hfill $\square$

\section{Proof of Theorem~\ref{main theorem}}\label{Proof of the main theorem}

\noindent
As in the previous section we assume that $\mcM$ {\bf \em is countable, binary, primitive, homogeneous and simple},
so by Fact~\ref{binary homogeneous simple structures are 1-based}, 
{\bf \em $\mcM$ is supersimple with finite SU-rank, 1-based and has trivial dependence.}
Just as in the previous section we  {\bf \em assume that the SU-rank of $\mcM$ is $\rho$.}
Furthermore, we assume that {\bf \em $\rho \geq 2$ and $C \subseteq M\meq$ is as in Assumption~\ref{assumption about C}}.
In addition, we adopt the following:

\begin{assump}\label{assumption about R}
For the rest of this section we assume that {\bf \em $\mcR$ is a binary random structure such that
(i)--(iv) of Lemma~\ref{getting suitable R} are satisfied}, so in particular $\mcC$ is a reduct of $\mcR$, 
which implies that $C = R$ (where $R$ is the universe of $\mcR$).
\end{assump}

\begin{rem}\label{remarks on R and C}{\rm
(i) A direct consequence of Assumption~\ref{assumption about R} is that if $\bar{c}, \bar{c}' \in C$ and
$\tp_\mcR(\bar{c}) = \tp_\mcR(\bar{c}')$ then $\tp_{\mcM\meq}(\bar{c}) = \tp_{\mcM\meq}(\bar{c}')$.\\
(ii) Let $p \in S^\mcR_1(\es)$ and let 
$X = \{c \in R : \mcR \models p(c)\}$.
By Theorem~\ref{characterization of equivalence relations in a special case},
it follows that there is no nontrivial equivalence relation on $X$ which is $\es$-definable in $\mcR$.
This can also be proved directly by a straightforward argument.
(The existence of a nontrivial $\es$-definable equivalence relation that is not definable by a unary formula
would contradict the defining property of a binary random structure.)
}\end{rem}

\noindent
Given this framework, including the facts of the previous section,
Theorem~\ref{main theorem} is a consequence of 
Proposition~\ref{there is a unique 1-type over acl(empty set)} (via its corollary)
and Lemma~\ref{coordinates are not in the definable closure}.
Most of the work is devoted to proving Proposition~\ref{there is a unique 1-type over acl(empty set)}.

\begin{lem}\label{coordinates are not in the definable closure}
(i) There is no formula $\varphi(x,y)$ (without parameters) such that for some $a \in M$ there are
$c, c' \in \crd(a)$ such that $\mcM\meq \models \varphi(c, a) \wedge \neg\varphi(c', a)$.\\
(ii) For every $a \in M$ and every $c \in\crd(a)$, $c \notin \dcl_{\mcM\meq}(a)$.\\
(iii) For every $a \in M$ and all $c, c' \in \crd(a)$ there are an ordering $c_1, \ldots, c_{\rho-1}$ of
$\crd(a) \setminus \{c\}$ and an ordering $c'_1, \ldots, c'_{\rho-1}$ of $\crd(a) \setminus \{c'\}$ such that
\[
(a, c, c_1, \ldots, c_{\rho-1}) \equiv_{\mcM\meq} (a, c', c'_1, \ldots, c'_{\rho-1}).
\]
\end{lem}

\noindent
{\bf Proof.}
(i) Suppose that there are a formula $\varphi(x, y)$, $a \in M$ and $c, c' \in \crd(a)$ such that
$\mcM\meq \models  \varphi(c, a) \wedge \neg\varphi(c', a)$.
Since $\mcM$ is primitive it follows that $\tp_\mcM(b) = \tp_\mcM(b')$ for all $b, b' \in M$ and therefore the previous statement
holds for {\em all} $a \in M$. Moreover, by primitivity it follows that for all $a, b \in M$, 
\[
|\{c \in \crd(a) : \mcM\meq \models \varphi(c, a)\}| \ = \
|\{c \in \crd(b) : \mcM\meq \models \varphi(c, b)\}|.
\]
This implies that the following is an equivalence relation on $M$, which is $\es$-definable in $\mcM$:
\[
x \sim y \ \Longleftrightarrow \ \forall z \Big( \big(z \in \crd(x) \ \wedge \ \varphi(z, x) \big)
\ \rightarrow \ 
\big( z \in \crd(y) \ \wedge \ \varphi(z, y) \big) \Big).
\]
Take $a, b \in M$ such that $a \ind b$.
By Fact~\ref{existence of C}~(ii)--(iv), $\crd(a) \cap \crd(b) = \es$, so $a \not\sim b$.
Hence `$\sim$' has at least two classes (actually infinitely many).
Let $\varphi(\mcM\meq, a) \cap C = \{c_1, \ldots, c_k\}$ and
$\crd(a) = \{c_1, \ldots, c_\rho\}$, where by assumption $\rho > k > 0$.
From Fact~\ref{existence of C} it follows that
$a, c_{k+1}, \ldots, c_\rho \notin \acl_{\mcM\meq}(c_1, \ldots, c_k)$, so there are 
(by the existence of nonforking extensions, for example)
$a', c'_{k+1}, \ldots, c'_\rho \notin \acl_{\mcM\meq}(a)$ such that
\[
(a, c_1, \ldots, c_\rho) \equiv_{\mcM\meq} (a', c_1, \ldots, c_k, c'_{k+1}, \ldots, c'_\rho).
\]
Then $a' \neq a$ and $a' \sim a$, so `$\sim$' is nontrivial.
This contradicts that $\mcM$ is primitive. Hence~(i) is proved.

Part~(ii) follows directly from~(i) because we assume that that the SU-rank $\rho$ is at least two 
and hence $|\crd(a)| = \rho \geq 2$ for every $a \in M$.

(iii) Let $a \in M$ and $c, c' \in \crd(a)$. Let $\varphi(x, y)$ isolate $\tp_{\mcM\meq}(c, a)$.
By~(i), $\mcM\meq \models \varphi(c', a)$, so $(a, c) \equiv_{\mcM\meq} (a, c')$.
Since $\crd(a)$ is finite and $\{a\}$-definable it follows that there are orderings 
$c_1, \ldots, c_{\rho-1}$ and $c'_1, \ldots, c'_{\rho-1}$ of 
$\crd(a) \setminus \{c\}$ and $\crd(a) \setminus \{c'\}$, respectively, such that
\[
(a, c, c_1, \ldots, c_{\rho-1}) \equiv_{\mcM\meq} (a, c', c'_1, \ldots, c'_{\rho-1}).
\]
\hfill $\square$
\\

\noindent
The main part of the proof of Theorem~\ref{main theorem} consists of proving the following:

\begin{prop}\label{there is a unique 1-type over acl(empty set)}
For all $c, c' \in C$, $\tp_{\mcM\meq}(c / \acl_{\mcM\meq}(\es)) = \tp_{\mcM\meq}(c' / \acl_{\mcM\meq}(\es))$.
\end{prop}

\noindent
We postpone the proof of Proposition~\ref{there is a unique 1-type over acl(empty set)} for a little while and
first show how it is used to prove Theorem~\ref{main theorem} via the following corollary.

\begin{cor}\label{corollary to there is a unique 1-type}
For all $0 < n < \omega$ and all $c_1, \ldots, c_n, c'_1, \ldots, c'_n \in C$,
\[
(c_1, \ldots, c_n) \equiv_{\mcM\meq} (c'_1, \ldots, c'_n) \ \text{ if and only if } \
(c_1, \ldots, c_n) \equiv_\mcR (c'_1, \ldots, c'_n).
\]
\end{cor}

\noindent
{\bf Proof.}
By Proposition~\ref{there is a unique 1-type over acl(empty set)},
there is only one $p \in S^{\mcM\meq}_1(\acl_{\mcM\meq}(\es))$ which is realized in $C$.
From part~(ii) of Fact~\ref{getting suitable R} it follows that the vocabulary of $\mcR$ has only one unary relation symbol $P$
and $\mcR \models \forall x P(x)$. 
Suppose that 
$(c_1, \ldots, c_n) \equiv_{\mcM\meq} (c'_1, \ldots, c'_n)$, so in particular
$(c_i, c_j) \equiv_{\mcM\meq} (c'_i, c'_j)$ for all $i$ and $j$.
By Fact~\ref{getting suitable R}~(iii) it follows that $(c_i, c_j) \equiv_\mcR (c'_i, c'_j)$ for all $i$ and $j$.
Since $\mcR$ is binary and has elimination of quantifiers we get 
$(c_1, \ldots, c_n) \equiv_\mcR (c'_1, \ldots, c'_n)$.
The other direction follows from the assumption that $\mcC$ is a reduct of $\mcR$
(and was stated in Remark~\ref{remarks on R and C}).
\hfill $\square$
\\

\noindent
We now show how Lemma~\ref{coordinates are not in the definable closure}
and
Corollary~\ref{corollary to there is a unique 1-type} 
imply our main result:

\medskip

\noindent
{\bf Theorem~\ref{main theorem}}
{\em 
Suppose that $\mcM$ is countable, binary, homogeneous, primitive and simple.
Then the SU-rank of $Th(\mcM)$ is 1.
}

\medskip

\noindent
{\bf Proof.}
By Fact~\ref{binary homogeneous simple structures are 1-based}, the premises of 
the theorem imply that
$\mcM$ satisfies all conditions assumed in this section.
Suppose, as in this whole section, that the SU-rank of $Th(\mcM)$ is $\rho \geq 2$.
Then, for every $a \in M$, $\su(a) = \rho$ and
therefore (by Fact~\ref{derived properties of C and crd}~(i)) $|\crd(a)| = \rho$. 
By Fact~\ref{lines are rigid in R} and
Corollary~\ref{corollary to there is a unique 1-type} we get:
\begin{itemize}
\item[] For every line $\{c_1, \ldots, c_\rho\} \subseteq C$ and every
nontrivial permutation $\pi$ of $\{1, \ldots, \rho\}$,
$(c_1, \ldots, c_\rho) \not\equiv_{\mcM\meq} (c_{\pi(1)}, \ldots, c_{\pi(\rho)})$.
\end{itemize}
This implies that for every $a \in M$ and every $c \in \crd(a)$, $c \in \dcl_{\mcM\meq}(a)$.
But this contradicts Lemma~\ref{coordinates are not in the definable closure}.
\hfill $\square$
\\

\noindent
The rest of this section is devoted to proving Proposition~\ref{there is a unique 1-type over acl(empty set)}.
We begin with a sequence of lemmas, 
numbered from~\ref{lines with intersecting the same E-classes are isomorphic}
to~\ref{lines intersect all E-classes}, 
which deal with properties of coordinates in the present context and of the equivalence
relation $\tp_{\mcM\meq}(x / \acl_{\mcM\meq}(\es)) = \tp_{\mcM\meq}(y / \acl_{\mcM\meq}(\es))$
restricted to $C$.
Then we have the tools to finish the proof of Proposition~\ref{there is a unique 1-type over acl(empty set)};
this part begins with Notation~\ref{notation for the rest of section showing that E has only one class}.

\noindent
For the rest of this section we use the following definition and notation:

\begin{defin}\label{definition of E}{\rm
For all $c, c' \in C$, 
\[
E(c, c') \ \Longleftrightarrow \ 
\tp_{\mcM\meq}(c / \acl_{\mcM\meq}(\es)) \ = \ \tp_{\mcM\meq}(c' / \acl_{\mcM\meq}(\es)).
\]
}\end{defin}

\noindent
Observe that, since  $\mcM$ is $\omega$-categorical, $E$ is $\es$-definable in $\mcM\meq$.

\begin{lem}\label{lines with intersecting the same E-classes are isomorphic}
Suppose that $c_1, \ldots, c_\rho, c'_1, \ldots, c'_\rho \in C$. 
If
\[
(c_1, \ldots, c_\rho) \equiv_{\mcM\meq} (c'_1, \ldots, c'_\rho) \
\text{ and \ $E(c_i, c'_i)$ for all $i = 1, \ldots, \rho$,}
\]
then $(c_1, \ldots, c_\rho) \equiv_\mcR (c'_1, \ldots, c'_\rho)$.
\end{lem}

\noindent
{\bf Proof.}
From part~(ii) of Fact~\ref{getting suitable R} we get
$c_i \equiv_\mcR c'_i$ for every $i = 1, \ldots, \rho$.
This together with parts~(iii) and~(iv) of  Fact~\ref{getting suitable R} gives
$(c_i, c_j) \equiv_\mcR (c'_i, c_j)$ for all $i = 1, \ldots, \rho$.
Since the vocabulary of $\mcR$ is binary and $\mcR$ has elimination of quantifiers we get 
$(c_1, \ldots, c_\rho) \equiv_\mcR (c'_1, \ldots, c'_\rho)$.
\hfill $\square$

\begin{lem}\label{lines and E-classes}
Let $A \subseteq C$ be a line and let $X, Y \subseteq C$ be $E$-classes.
If $A \cap X \neq \es$ and $A \cap Y \neq \es$, then $|A \cap X| = |A \cap Y|$.
\end{lem}

\noindent
{\bf Proof.}
Suppose for a contradiction that there are a line $A$ and $E$-classes $X, Y$ such that 
$A \cap X \neq \es$, $A \cap Y \neq \es$ and $k = |A \cap X| \neq |A \cap Y|$.
Let $a \in M$ be such that $\crd(a) = A$, let $c \in A \cap X$ and $c' \in A \cap Y$.
Let $\varphi(x, y)$ be the formula which expresses that
\begin{itemize}
\item[] ``$x \in \crd(y)$ and there are exactly $k$ different elements $z \in \crd(y)$ such that $E(x, z)$''.
\end{itemize}
Then $\mcM\meq \models \varphi(c, a) \wedge \neg\varphi(c', a)$, which contradicts 
Lemma~\ref{coordinates are not in the definable closure}~(i).
\hfill $\square$

\begin{lem}\label{lines intersect all E-classes}
There is a  number $s \geq 1$ such that for every 
line $A \subseteq C$ and every $E$-class $X \subseteq C$, $|A \cap X| = s$.
\end{lem}

\noindent
{\bf Proof.}
By Lemma~\ref{lines and E-classes} and since all elements of $M$ have the same type
it suffices to show that every line has nonempty intersection with every $E$-class.
Let $l$ be the number of $E$-classes (so $l < \omega$).
Let $A$ be any line and let $k$ be the number of $E$-classes $X$ such that $A \cap X \neq \es$.
Since all elements of $\mcM$ realize the same 1-type over $\es$ it follows that for every line $A'$ there are exactly $k$
$E$-classes $X$ such that $A' \cap X \neq \es$.

Consider then following $\es$-definable (in $\mcM$) 
equivalence relation on $M$ :
\begin{align*}
x \sim y \ \Longleftrightarrow \ 
&\text{ for all $u \in \crd(x)$ and all $v \in \crd(y)$ there are}\\
&\text{ $u' \in \crd(y)$ and $v' \in \crd(x)$
such that $E(u, u')$ and $E(v, v')$.}
\end{align*}
We will show that if $l > k$ then `$\sim$' is nontrivial, contradicting that $\mcM$ is primitive.
So suppose that $l > k$. 
Take any $a \in M$ and let $\crd(a) = \{c_1, \ldots, c_\rho\}$.
As $l > k$ there is $c' \in C$ such that $\neg E(c', c_i)$ for all $i = 1, \ldots, \rho$.
By Fact~\ref{existence of C}~(v) there is $a' \in M$ such that $c' \in \crd(a')$.
Then $a \not\sim a'$ so `$\sim$' has at least two classes.
Let $X$ be the $E$-class of $c_1$.
As $c_1 \notin \acl_{\mcM\meq}(c_2, \ldots, c_\rho)$ there is $c'_1 \in X \setminus \acl_{\mcM\meq}(c_1, \ldots, c_\rho)$
such that
\[
(c'_1, c_2, \ldots, c_\rho) \equiv_{\mcM\meq} (c_1, c_2, \ldots, c_\rho).
\]
Then $\{c'_1, c_2, \ldots, c_\rho\}$ is a line so $\{c'_1, c_2, \ldots, c_\rho\} = \crd(a')$ for some $a' \in M$.
It follows that $a' \neq a$ and $a' \sim a$. Thus `$\sim$' is nontrivial, a contradiction.

Hence we conclude that $l = k$ which implies that for every line
$A \subseteq C$ and every $E$-class $X \subseteq C$, $A \cap X \neq \es$.
\hfill $\square$

\begin{notation}\label{notation for the rest of section showing that E has only one class}{\rm
For the rest of this section we will use the following notation:
\begin{itemize}
\item[(i)] Let $X_1, \ldots, X_l$ enumerate all $E$-classes (without repetition).

\item[(ii)] Let $A_0$ be any line.
By Lemma~\ref{lines intersect all E-classes}, there is $s \geq 1$ such that 
$|A_0 \cap X_i| = s$ for every $i = 1, \ldots, l$.
For $i = 1, \ldots, l$, let $\bar{d}_i = (d_{i, 1}, \ldots, d_{i,s})$ enumerate $A_0 \cap X_i$.
By Lemma~\ref{coordinates are not in the definable closure} 
we may assume that $\bar{d}_i \equiv_{\mcM\meq} \bar{d}_j$ for all $i$ and $j$.

\item[(iii)] Let
\begin{align*}
p(\bar{x}_1 \ldots \bar{x}_l) &= \tp_{\mcM\meq}(\bar{d}_1 \ldots \bar{d}_l) \ \text{ and}\\
p^+(\bar{x}_1 \ldots \bar{x}_l) &= \tp_\mcR(\bar{d}_1 \ldots \bar{d}_l).
\end{align*}
\end{itemize}
}\end{notation}

\noindent
Observe that, by Remark~\ref{remarks on R and C}~(i),
for all $\bar{c}_1, \ldots, \bar{c}_l \in C$, if $\mcR \models p^+(\bar{c}_1 \ldots \bar{c}_l)$
then $\mcM\meq \models p(\bar{c}_1 \ldots \bar{c}_l)$.
Morover, 
by Remark~\ref{remarks on R and C}~(ii),
if $\mcR \models p^+(\bar{c}_1 \ldots \bar{c}_l)$ then $\rng(\bar{c}_i) \subseteq X_i$ for all $i = 1, \ldots, l$,
and if $\mcM\meq \models p(\bar{c}_1 \ldots \bar{c}_l)$ then, for every $i = 1, \ldots, l$, all members of $\bar{c}_i$ belong
to the same $E$-class.

\begin{lem}\label{all lines can be ordered to realize p}
Every line $A \subseteq C$ can be enumerated as 
$\bar{c}_1 \ldots \bar{c}_l$, where $\bar{c}_i = (c_{i,1}, \ldots c_{i,s})$ for $i = 1, \ldots, l$, so that
$\mcM\meq \models p(\bar{c}_1 \ldots \bar{c}_l)$ and, for all $i = 1, \ldots, l$ and all $j = 1, \ldots, s$,
$c_{i,j} \in X_i$.
\end{lem}

\noindent
{\bf Proof.}
Define on $M$:
\begin{align*}
a \sim a' \ \Longleftrightarrow \
&\text{ $\crd(a)$ and $\crd(a')$ can be enumerated as 
$\bar{c}_1 \ldots \bar{c}_l$ and $\bar{c}'_1 \ldots \bar{c}'_l$,} \\
&\text{ respectively, so that
$\mcM\meq \models p(\bar{c}_1 \ldots \bar{c}_l) \wedge p(\bar{c}'_1 \ldots \bar{c}'_l)$ and,} \\
&\text{ for all $i = 1, \ldots, l$ and all $j = 1, \ldots, s$, $E(c_{i, j}, c'_{i, j})$.}
\end{align*}
Then `$\sim$' is clearly $\es$-definable in $\mcM$
(because by the $\omega$-categoricity of $\mcM$, $p$ is isolated), 
as well as reflexive and symmetric.
We will show that `$\sim$' is transitive, hence an equivalence relation.
Since $\mcR$ is a binary random structure and $\mcC$ is a reduct of $\mcR$ it is easy to see that
there are distinct $a, a' \in M$ such that $a \sim a'$.
As $\mcM$ is primitive it follows that $a \sim a'$ for all $a, a' \in M$.
The conclusion of the lemma follows from this. Hence it remains to show that `$\sim$' is transitive.

Suppose that $a \sim a' \sim a''$.
From $a \sim a'$ it follows that $\crd(a)$ and $\crd(a')$ can be enumerated as
$\bar{c}_1 \ldots \bar{c}_l$ and $\bar{c}'_1 \ldots \bar{c}'_l$, respectively, so that
\begin{align}\label{relationship between a and a'}
&\text{ $\mcM\meq \models p(\bar{c}_1 \ldots \bar{c}_l) \wedge p(\bar{c}'_1 \ldots \bar{c}'_l)$ and,} \\
&\text{ for all $i = 1, \ldots, l$ and all $j = 1, \ldots, s$, $E(c_{i, j}, c'_{i, j})$.} \nonumber
\end{align}
From $a' \sim a''$ it follows that $\crd(a')$ and $\crd(a'')$ can be enumerated as
$\bar{c}^*_1 \ldots \bar{c}^*_l$ and $\bar{c}''_1 \ldots \bar{c}''_l$, respectively, so that
\begin{align}\label{relationship between a' and a''}
&\text{ $\mcM\meq \models p(\bar{c}^*_1 \ldots \bar{c}^*_l) \wedge p(\bar{c}''_1 \ldots \bar{c}''_l)$ and,} \\
&\text{ for all $i = 1, \ldots, l$ and all $j = 1, \ldots, s$, $E(c^*_{i, j}, c''_{i, j})$.} \nonumber
\end{align}

Since $\bigcup_{i=1}^l \rng(\bar{c}'_i) = \bigcup_{i=1}^l \rng(\bar{c}^*_i)$ it follows 
from~(\ref{relationship between a and a'}) and~(\ref{relationship between a' and a''})
that there is a 
\[
\text{permutation $\pi$ of } \big\{ (i, j) : 1 \leq i \leq l, \ 1 \leq j \leq s \big\}
\]
such that 
\begin{equation}\label{pi-c' and c* are equivalent}
\text{ for all $i = 1, \ldots, l$ and all $j = 1, \ldots, s$, $c'_{\pi(i,j)} = c^*_{i,j}$}
\end{equation}
and (since $p$ is a complete type)
\begin{align*}
\mcM\meq \models \forall x_{1, 1}, \ldots, x_{1, s}, \ldots, x_{l, 1}, \ldots, x_{l, s} \Big(
p(x_{1, 1}, \ldots, x_{1, s}, \ldots, x_{l, 1}, \ldots, x_{l, s}) \ \rightarrow \\
p(x_{\pi(1, 1)}, \ldots, x_{\pi(1, s)}, \ldots, x_{\pi(l, 1)}, \ldots, x_{\pi(l, s)}) \Big).
\end{align*}
This together with~(\ref{relationship between a and a'}) and~(\ref{relationship between a' and a''}) gives
\begin{align}\label{pi-c and c''}
\mcM\meq \models \
&p(c_{\pi(1, 1)}, \ldots, c_{\pi(1, s)}, \ldots, c_{\pi(l, 1)}, \ldots, c_{\pi(l, s)}) \\ 
\wedge \
&p(c''_{1, 1}, \ldots, c''_{1, s}, \ldots, c''_{l, 1}, \ldots, c''_{l, s}). \nonumber
\end{align}
From~(\ref{relationship between a' and a''}) and~(\ref{pi-c' and c* are equivalent})
it follows that, for all $i = 1, \ldots, l$ and all $j = 1, \ldots, s$, 
$E(c'_{\pi(i,j)}, c''_{i, j})$.
From~(\ref{relationship between a and a'}) we get 
$E(c_{i, j}, c'_{i, j})$ and hence
$E(c_{\pi(i,j)}, c'_{\pi(i,j)})$ for all
$i = 1, \ldots, l$ and all $j = 1, \ldots, s$.
By transitivity of $E$ we get 
$E(c_{\pi(i,j)}, c''_{i, j})$ for all $i = 1, \ldots, l$ and all $j = 1, \ldots, s$.
This and~(\ref{pi-c and c''}) imply that $a \sim a''$, so `$\sim$' is transitive.
\hfill $\square$
\\

\noindent
Recall that if $\mcR \models p^+(\bar{c}_1 \ldots \bar{c}_l)$ then
$\mcM\meq \models p(\bar{c}_1 \ldots \bar{c}_l)$  and, for all $i = 1, \ldots, l$
and all $j = 1, \ldots, s$, $c_{i,j} \in X_i$. Moreover, the range of every tuple that realizes $p$ is a line.
\\

\noindent
{\bf \em For the rest of this section we assume that $l > 1$.} 
\\

\noindent
From this we will derive a contradiction and thus prove Proposition~\ref{there is a unique 1-type over acl(empty set)}.
Since $\mcR$ is a binary random structure there are
\begin{align}\label{properties of the a-i and b-i}
&\bar{a}_i = (a_{i,1}, \ldots, a_{i,s}), i = 1, \ldots, l, \\
&\bar{b}_i = (b_{i,1}, \ldots, b_{i,s}), i = 1, \ldots, l, \ \text{  such that} \nonumber \\
&\bar{a}_l = \bar{b}_l, \ 
\bigcup_{i=1}^l \rng(\bar{a}_i) \ \cap \ \bigcup_{i=1}^l \rng(\bar{b}_i) = \rng(\bar{a}_l) \ \text{ and} \nonumber \\
&\mcR \models p^+(\bar{a}_1 \ldots \bar{a}_l) \ \wedge \ p^+(\bar{b}_1 \ldots \bar{b}_l). \nonumber
\end{align}

\noindent
Moreover, there is a disjoint copy (up to isomorphism in $\mcR$) of the above elements.
More precisely, there are
\begin{align}\label{properties of the a'-i and b'-i}
&\bar{a}'_i = (a'_{i,1}, \ldots, a'_{i,s}), i = 1, \ldots, l, \\
&\bar{b}'_i = (b'_{i,1}, \ldots, b'_{i,s}), i = 1, \ldots, l, \ \text{  such that} \nonumber \\
&(\bar{a}'_1, \ldots, \bar{a}'_l, \bar{b}'_1, \ldots, \bar{b}'_l) \equiv_\mcR
(\bar{a}_1, \ldots, \bar{a}_l, \bar{b}_1, \ldots, \bar{b}_l), \ \text{ and} \nonumber \\
&\bigcup_{i=1}^l (\rng(\bar{a}'_i) \cup \rng(\bar{b}'_i)) \ \cap \ 
\bigcup_{i=1}^l (\rng(\bar{a}_i) \cup \rng(\bar{b}_i)) = \es, \nonumber
\end{align}
and consequently
\begin{align*}
&\bar{a}'_l = \bar{b}'_l, \ 
\bigcup_{i=1}^l \rng(\bar{a}'_i) \ \cap \ \bigcup_{i=1}^l \rng(\bar{b}'_i) = \rng(\bar{a}'_l) \ \text{ and} \\
&\mcR \models p^+(\bar{a}'_1 \ldots \bar{a}'_l) \ \wedge \ p^+(\bar{b}'_1 \ldots \bar{b}'_l). \nonumber
\end{align*}

\noindent
Moreover (as $\mcR$ is a binary random structure),
we can choose $\bar{a}'_i, \bar{b}'_i$, $i = 1, \ldots,l$, so that, in addition
to~(\ref{properties of the a'-i and b'-i}), 
the following holds:
\begin{align}\label{a'-1 and b'-1 have the same type over a-i and b-i}
\text{for all } i = 1, \ldots, l, \ 
(\bar{a}'_1, \bar{a}_i) \equiv_\mcR (\bar{b}'_1, \bar{a}_i) \ \text{ and } \
(\bar{a}'_1, \bar{b}_i) \equiv_\mcR (\bar{b}'_1, \bar{b}_i).
\end{align}

\noindent
If $\sigma$ is a permutation of a set $I$, $\{e_i : i \in I\}$ is a set indexed by $I$
and $\bar{e} = (e_{i_1}, \ldots, e_{i_k})$ where 
$i_1, \ldots, i_k \in I$, then we let $\sigma(\bar{e})$ denote the sequence $(e_{\sigma(i_1)}, \ldots, e_{\sigma(i_k)})$.

\begin{lem}\label{the existence of a''-i and b''-i}
There are 
\begin{align}\label{properties of the a''-i and b''-i}
&\bar{a}''_i = (a''_{i,1}, \ldots, a''_{i,s}), i = 1, \ldots, l, \\
&\bar{b}''_i = (b''_{i,1}, \ldots, b''_{i,s}), i = 1, \ldots, l, \ \text{ and} \nonumber \\
&\text{a permutation $\sigma$ of $\{(i, j) : 1 \leq i \leq l, \ 1 \leq j \leq s\}$  such that} \nonumber \\
&\bar{a}''_1 = \bar{b}''_1, \ \bigcup_{i=1}^l \rng(\bar{a}''_i) \ \cap \ 
\bigcup_{i=1}^l \rng(\bar{b}''_i) = \rng(\bar{a}''_1), \nonumber \\
&\mcR \models p^+(\bar{a}''_1 \ldots \bar{a}''_l) \ \wedge \ p^+(\bar{b}''_1 \ldots \bar{b}''_l) \ \text{ and} \nonumber \\
&(\bar{a}_1, \ldots, \bar{a}_l, \bar{b}_1, \ldots, \bar{b}_l) \equiv_{\mcM\meq}
(\sigma(\bar{a}''_1), \ldots, \sigma(\bar{a}''_l), \sigma(\bar{b}''_1), \ldots, \sigma(\bar{b}''_l)). \nonumber
\end{align}
\end{lem}

\noindent
{\bf Proof.}
From the choice of $p$ in Notation~\ref{notation for the rest of section showing that E has only one class}, 
it follows that $\bar{a}_1 \equiv_{\mcM\meq} \bar{a}_l$.
Hence there is an automorphism $f$ of $\mcM\meq$ such that $f(\bar{a}_l) = \bar{a}_1$
(and since $\bar{a}_l = \bar{b}_l$ we also have $f(\bar{b}_l) = \bar{a}_1$).
Let 
\begin{align*}
&\bar{a}^*_i = (a^*_{i,1}, \ldots, a^*_{i,s}) = (f(a_{i,1}), \ldots, f(a_{i,s})) \ \text{ and} \\
&\bar{b}^*_i = (b^*_{i,1}, \ldots, b^*_{i,s}) = (f(b_{i,1}), \ldots, f(b_{i,s})) \ \text{ for } i = 1, \ldots, l.
\end{align*}
Then
\begin{equation}\label{equivalence in Meq between a-i b-i and a*-i b*-i}
(\bar{a}_1, \ldots, \bar{a}_l, \bar{b}_1, \ldots, \bar{b}_l) \ \equiv_{\mcM\meq} \
(\bar{a}^*_1, \ldots, \bar{a}^*_l, \bar{b}^*_1, \ldots, \bar{b}^*_l) \ 
\text{ where $\bar{a}^*_l = \bar{b}^*_l = \bar{a}_1$},
\end{equation}
from which it follows that
\begin{align}\label{a*-i and b*-i realize p}
&\mcM\meq \models p(\bar{a}^*_1 \ldots \bar{a}^*_l) \ \wedge \ p(\bar{b}^*_1 \ldots \bar{b}^*_l) \ \text{ and} \\
&\text{for all $i = 1, \ldots, l$ and all $j = 1, \ldots, s$, } E(a^*_{i,j}, b^*_{i,j}). \nonumber
\end{align}
Then $\rng(\bar{a}^*_1) \cup \ldots \cup \rng(\bar{a}^*_l)$ is a line so, 
by Lemma~\ref{all lines can be ordered to realize p},
there is a 
\[
\text{permutation
$\pi$ of $\{(i,j) : 1 \leq i \leq l, \ 1 \leq j \leq s\}$}
\]
such that
\begin{align}\label{pi-a*-i satisfy p}
&\text{$\rng(\pi(\bar{a}^*_i)) \subseteq X_i$ for all $i = 1, \ldots, l$, and} \\
&\mcM\meq \models p(\pi(\bar{a}^*_1) \ldots \pi(\bar{a}^*_l)). \nonumber
\end{align}

\noindent
This and~(\ref{a*-i and b*-i realize p}) implies that 
\begin{equation}\label{pi-b*-i realizes p}
\mcM\meq \models p(\pi(\bar{b}^*_1) \ldots \pi(\bar{b}^*_l)).
\end{equation}
From~(\ref{a*-i and b*-i realize p})
and~(\ref{pi-a*-i satisfy p})
it follows that
\begin{equation}\label{pi-b*-i is in X-i}
\text{for all $i = 1, \ldots, l$, $\rng(\pi(\bar{b}^*_i)) \subseteq X_i$.}
\end{equation}
Thus we have 
\[
\big\{ \pi(1, j) : 1 \leq j \leq l \big\} \ = \ \big\{ (l, j) : 1 \leq j \leq s \big\},
\]
so there is a permutation $\gamma$ of $\big\{(l, j) : 1 \leq i \leq s \big\}$ such that
\[
\bar{a}^*_1 = \pi^{-1}\gamma(\bar{a}^*_l) \ \text{ and } \
\bar{b}^*_1 = \pi^{-1}\gamma(\bar{b}^*_l).
\]
As $\bar{a}^*_l = \bar{b}^*_l$ we get $\bar{a}^*_1 = \bar{b}^*_1$ and hence 
\[
\pi(\bar{a}^*_1) = \pi(\bar{b}^*_1).
\]

\noindent
By~(\ref{pi-a*-i satisfy p}) --~(\ref{pi-b*-i is in X-i})
and
Lemma~\ref{lines with intersecting the same E-classes are isomorphic} we get
\begin{equation*}
\mcR \models p^+(\pi(\bar{a}^*_1) \ldots \pi(\bar{a}^*_l)) \ \wedge \
p^+(\pi(\bar{b}^*_1) \ldots \pi(\bar{b}^*_l)).
\end{equation*}
If we now let $\bar{a}''_i = \pi(\bar{a}^*_i)$, $\bar{b}''_i = \pi(\bar{b}^*_i)$, for $i = 1, \ldots, l$,
and $\sigma = \pi^{-1}$, then~(\ref{properties of the a''-i and b''-i}) is satisfied, so the lemma is proved.
\hfill $\square$
\\

\noindent
By Lemma~\ref{the existence of a''-i and b''-i} there are
$\bar{a}''_i, \bar{b}''_i$, $i = 1, \ldots, l$, so that~(\ref{properties of the a''-i and b''-i}) holds.
As $\mcR$ is a binary random structure, 
and by~(\ref{a'-1 and b'-1 have the same type over a-i and b-i}), 
we can choose these elements so that, in addition to~(\ref{properties of the a''-i and b''-i}),
\begin{align}\label{properties of a''-i and b''-i in relation to a-i b-i etc}
&\bigcup_{i=1}^l \big(\rng(\bar{a}''_i) \cup \rng(\bar{b}''_i)\big) \ \cap \ 
\bigcup_{i=1}^l \big(\rng(\bar{a}_i) \cup \rng(\bar{b}_i) \cup \rng(\bar{a}'_i) \cup \rng(\bar{b}'_i)\big)
= \es, \\
&(\bar{a}''_1, \ldots, \bar{a}''_l, \bar{a}_1, \ldots, \bar{a}_l) \ \equiv_\mcR \ 
(\bar{a}'_1, \ldots, \bar{a}'_l, \bar{a}_1, \ldots, \bar{a}_l), \nonumber \\
&(\bar{a}''_1, \ldots, \bar{a}''_l, \bar{b}_1, \ldots, \bar{b}_l) \ \equiv_\mcR \ 
(\bar{a}'_1, \ldots, \bar{a}'_l, \bar{b}_1, \ldots, \bar{b}_l), \nonumber \\
&(\bar{b}''_1, \ldots, \bar{b}''_l, \bar{a}_1, \ldots, \bar{a}_l) \ \equiv_\mcR \ 
(\bar{b}'_1, \ldots, \bar{b}'_l, \bar{a}_1, \ldots, \bar{a}_l), \ \text{ and} \nonumber \\
&(\bar{b}''_1, \ldots, \bar{b}''_l, \bar{b}_1, \ldots, \bar{b}_l) \ \equiv_\mcR \ 
(\bar{b}'_1, \ldots, \bar{b}'_l, \bar{b}_1, \ldots, \bar{b}_l). \nonumber
\end{align}

\noindent
Since $\mcC$ is a reduct of $\mcR$ (and $\mcC$ is canonically embedded in $\mcM\meq$) it follows that
\begin{equation}\label{properties of a''-i and b''-i in relation to a-i b-i etc holds in Meq}
\text{in (\ref{properties of a''-i and b''-i in relation to a-i b-i etc}) we can replace `$\equiv_\mcR$' by `$\equiv_{\mcM\meq}$'.}
\end{equation}

\noindent
Let $a, a', a'', b, b', b'' \in M$ be such that
\begin{align*}
\crd(a) &= \rng(\bar{a}_1) \cup \ldots \cup \rng(\bar{a}_l), 
&\crd(b) &= \rng(\bar{b}_1) \cup \ldots \cup \rng(\bar{b}_l), \\
\crd(a') &= \rng(\bar{a}'_1) \cup \ldots \cup \rng(\bar{a}'_l), 
&\crd(b') &= \rng(\bar{b}'_1) \cup \ldots \cup \rng(\bar{b}'_l), \\
\crd(a'') &= \rng(\bar{a}''_1) \cup \ldots \cup \rng(\bar{a}''_l) \ \text{ and }
&\crd(b'') &= \rng(\bar{b}''_1) \cup \ldots \cup \rng(\bar{b}''_l).
\end{align*}
We have $a \in \dcl_{\mcM\meq}(\crd(a))$ and similarly for $a', a'', b, b'$ and $b''$.
Therefore it follows 
from~(\ref{properties of a''-i and b''-i in relation to a-i b-i etc holds in Meq}) 
and~(\ref{properties of a''-i and b''-i in relation to a-i b-i etc}) that
\begin{align*}
(a', b') &\equiv_{\mcM\meq} (a'', b''), &(a', a) &\equiv_{\mcM\meq} (a'', a), \\
(a', b) &\equiv_{\mcM\meq} (a'', b), &(b', a) &\equiv_{\mcM\meq} (b'', a), \ \text{ and} \\
(b', b) &\equiv_{\mcM\meq} (b'', b). &  &
\end{align*}

\noindent
Since all the involved elements belong to $M$ we can replace `$\equiv_{\mcM\meq}$' by `$\equiv_\mcM$'.
As $\mcM$ is binary with elimination of quantifiers we get
\begin{equation}\label{aba'b' is elementarily equivalent to aba''b''}
(a, b, a', b') \equiv_\mcM (a, b, a'', b'').
\end{equation}
Now consider a formula $\varphi(x_1, x_2, x_3, x_4)$ in the language of $\mcM\meq$ which expresses the following:
\begin{itemize}
\item[] ``$x_1, x_2, x_3, x_4 \in M$ and there is an $E$-class $X$ such that 
$\crd(x_1) \cap X = \crd(x_2) \cap X$ and $\crd(x_3) \cap X = \crd(x_4) \cap X$.''
\end{itemize}
It is straightforward to verify that 
\[
\mcM\meq \models \varphi(a, b, a', b') \ \wedge \ \neg\varphi(a, b, a'', b''),
\]
so $(a, b, a', b') \not\equiv_{\mcM\meq} (a, b, a'', b'')$. Since $a, a', a'', b, b', b'' \in M$ it follows that
$(a, b, a', b') \not\equiv_\mcM (a, b, a'', b'')$, which contradicts~(\ref{aba'b' is elementarily equivalent to aba''b''}).
This completes the proof of Proposition~\ref{there is a unique 1-type over acl(empty set)}
and hence of Theorem~\ref{main theorem}.

\section{Definable equivalence relations}\label{Definable equivalence relations}

\noindent
In this section we prove a result, Theorem~\ref{characterization of equivalence relations for omega-categorical},
about $\es$-definable equivalence relations on $n$-tuples (for any fixed $0 < n < \omega$)
in $\omega$-categorical supersimple structures with SU-rank 1 and
degenerate algebraic closure. 
One reason for doing this is that the author thinks that this result may be useful in future research about
nonbinary simple homogeneous structures.
Another reason is that a variant of Theorem~\ref{characterization of equivalence relations for omega-categorical},
namely Theorem~\ref{characterization of equivalence relations in a special case}, gives (under certain conditions)
a full characterization of the $\es$-definable equivalence relations on $n$-tuples,
for any $0 < n < \omega$. Theorem~\ref{characterization of equivalence relations in a special case}
is then used to prove 
Proposition~\ref{the generic tetrahedron-free hypergraph is 1-based}
which implies that the generic tetrahedron-free 3-hypergraph is 1-based.

{\bf \em Throughout this section we suppose that 
$\mcM$ is $\omega$-categorical, supersimple with SU-rank 1 and with degenerate algebraic closure.}
Let $n < \omega$ and let $p(\bar{x}) \in S^{\mcM\meq}_n(\acl_{\mcM\meq}(\es))$ be 
realized by some $n$-tuple of elements from $M$.

Let
\[
X = \{\bar{a} \in M^n : \mcM \models p(\bar{a}) \}.
\]
{\bf \em Suppose that $E$ is an equivalence relation on $X$ which is $\es$-definable in $\mcM$.}
In other words, there is a formula $\varphi(\bar{x}, \bar{y})$ 
without parameters, in the language of $\mcM$,
such that for all $\bar{a}, \bar{b} \in X$, $E(\bar{a}, \bar{b})$ if and only if $\mcM \models \varphi(\bar{a}, \bar{b})$.

\begin{theor}\label{characterization of equivalence relations for omega-categorical}
Suppose that $E$ is nontrivial, i.e. it has at least two equivalence classes and at least one of the equivalence classes
has more than one element.
Then there is a nonempty $I \subseteq \{1, \ldots, n\}$,
a group of permutations $\Gamma$ of $I$ and
$\es$-definable equivalence relations $E'$ and $E''$ on $X$ such that the following hold:
\begin{itemize}
\item[(a)] $E'' \subseteq E \subseteq E'$.

\item[(b)] For all $\bar{a} = (a_1, \ldots, a_n), \bar{b} = (b_1, \ldots, b_n) \in X$, $E'(\bar{a}, \bar{b})$ if and only if
there is a permutation $\gamma \in \Gamma$ of $I$
such that $a_i = b_{\gamma(i)}$ for all $i \in I$.

\item[(c)] For all $\bar{a} = (a_1, \ldots, a_n), \bar{b} = (b_1, \ldots, b_n) \in X$, $E''(\bar{a}, \bar{b})$ if and only if
$a_i = b_i$ for all $i \in I$ and if $\{1, \ldots, n\} \setminus I = \{i_1, \ldots, i_m\}$, then
\[
\tp_{\mcM\meq}(a_{i_1}, \ldots, a_{i_m} /\acl_{\mcM\meq}(\{a_i : i \in I\})) =
\tp_{\mcM\meq}(b_{i_1}, \ldots, b_{i_m} /\acl_{\mcM\meq}(\{a_i : i \in I\})).
\]
\end{itemize}
\end{theor}

\noindent
The rest of this section proves this theorem.
Without loss of generality we assume that $p(\bar{x})$ implies $x_i \neq x_j$ for all $1 \leq i < j \leq n$.
In this section and the next we frequently abuse notation by notationally identifying `$\bar{a}$' and `$\rng(\bar{a})$'.
{\bf \em From now on suppose that $E$ is nontrivial.}
Several times in this section and the next we will use the following:

\begin{observation}\label{facts about independence in degenerate structures}
(i) For all $\bar{a}, \bar{b}, \bar{c} \in M$, $\bar{a} \underset{\bar{c}}{\ind} \bar{b}$ if and only if
$\bar{a} \cap \bar{b} \subseteq \bar{c}$.\\
(ii) For all $\bar{a}, \bar{b}, \bar{c} \in M$ there is $\bar{a}' \in M$ such that 
$\tp_\mcM(\bar{a}' / \acl_{\mcM\meq}(\bar{c})) = \tp_{\mcM\meq}(\bar{a} / \acl_{\mcM\meq}(\bar{c}))$
and $\bar{a}' \cap \bar{b} \subseteq \bar{c}$.
\end{observation}

\noindent
Both parts of Observation~\ref{facts about independence in degenerate structures}
are straightforward to show and hold under the assumptions in this section.
We note however that~(ii) is a direct consequence
of \cite[Proposition 1.5 (1)]{Adl} (which shows that `algebraic independence' satisfies `full existence').
Part~(i) is {\em not} used in the proof of Theorem~\ref{characterization of equivalence relations in a special case},
but it {\em is} used in the proof of Proposition~\ref{the generic tetrahedron-free hypergraph is 1-based}, which has even 
stronger assumptions than the present section. 
Part~(ii) holds under the assumptions of Theorem~\ref{characterization of equivalence relations in a special case}
and the assumptions of Proposition~\ref{the generic tetrahedron-free hypergraph is 1-based}.

\begin{lem}\label{first lemma on disjoint tuples for omega-categorical}
Either $E(\bar{a}, \bar{b})$ holds for all disjoint $\bar{a}, \bar{b} \in X$,
or $\neg E(\bar{a}, \bar{b})$ holds for all disjoint $\bar{a}, \bar{b} \in X$.
\end{lem}

\noindent
{\bf Proof.}
For a contradiction suppose that $\bar{a}, \bar{b} \in X$ are disjoint and $E(\bar{a}, \bar{b})$ and that
$\bar{c}, \bar{d} \in X$ are disjoint and $\neg E(\bar{c}, \bar{d})$.
As $\mcM$ is $\omega$-categorical it follows from the definition of $X$ that there is $\bar{b}' \in X$ such that
\[
\tp_\mcM(\bar{b}, \bar{b}') = \tp_\mcM(\bar{c}, \bar{d}).
\]
Then $\bar{b} \cap \bar{b}' = \es$.
As $\acl_\mcM$ is degenerate and $\mcM$ has SU-rank 1 we get $\bar{b} \ind \bar{b}'$. 
By assumption, $\bar{a}$ and $\bar{b}$ are disjoint, 
so (by Observation~\ref{facts about independence in degenerate structures}~(i)) $\bar{a} \ind \bar{b}$. By definition of $X$, 
$\tp_{\mcM\meq}(\bar{b} / \acl_{\mcM\meq}(\es)) = p = \tp_{\mcM\meq}(\bar{b}' / \acl_{\mcM\meq}(\es))$.
Therefore the independence theorem implies that there is $\bar{e}$ such that $\bar{e}$ realizes $p$, 
and hence $\bar{e} \in X$,
$\tp_\mcM(\bar{a}, \bar{e}) = \tp_\mcM(\bar{a}, \bar{b})$ and $\tp_\mcM(\bar{b}, \bar{e}) = \tp_\mcM(\bar{b}, \bar{b}')$.
This implies that $E(\bar{a}, \bar{e})$ and $\neg E(\bar{b}, \bar{e})$.
Together with the assumption that $E(\bar{a}, \bar{b})$ we have a contradiction to the symmetry and transitivity
of $E$.
\hfill $\square$

\begin{lem}\label{second lemma on disjoint tuples for omega-categorical}
For all disjoint $\bar{a}, \bar{b} \in X$ we have $\neg E(\bar{a}, \bar{b})$.
\end{lem}

\noindent
{\bf Proof.}
Suppose that the lemma is false. By 
Lemma~\ref{first lemma on disjoint tuples for omega-categorical},
for all disjoint $\bar{a}, \bar{b} \in X$ we have $E(\bar{a}, \bar{b})$.
Since $E$ is nontrivial there are $\bar{b}, \bar{c} \in X$ such that $\neg E(\bar{b}, \bar{c})$. 
By Observation~\ref{facts about independence in degenerate structures}~(ii), 
there is $\bar{a} \in X$ which is disjoint from $\bar{b}$ and from $\bar{c}$.
Then $E(\bar{a}, \bar{b})$, $E(\bar{a}, \bar{c})$ and $\neg E(\bar{b}, \bar{c})$, which contradicts the
symmetry and transitivity of $E$.
\hfill $\square$

\begin{lem}\label{the most technical lemma about intersection of tuples for omega-categorical}
Suppose that $\bar{a}, \bar{b} \in X$, $E(\bar{a}, \bar{b})$,  $\bar{a} \cap \bar{b} \neq \es$,
$\bar{a} \cap \bar{b} =
\{a_{i_1}, \ldots, a_{i_k}\} = \{b_{j_1}, \ldots, b_{j_k}\}$ (where elements are listed without repetition) and
$\{i_1, \ldots, i_k\} \neq \{j_1, \ldots, j_k\}$ (so $k \geq 1$).
Then there is $\bar{c} \in X$ such that $E(\bar{a}, \bar{c})$ and 
$\bar{a} \cap \bar{c}$ is a proper subset of $\bar{a} \cap \bar{b}$.
\end{lem}

\noindent
{\bf Proof.}
By reindexing variables if necessary we may, without loss of generality, assume that
for some $0 < m \leq n-k$ and some permutation $\gamma$ of $\{1, \ldots, k\}$
\[
\bar{a} \cap \bar{b} = \{a_1, \ldots, a_k\}  \ 
\text{ and $a_i = b_{m+\gamma(i)}$ for all $i = 1, \ldots ,k$.}
\]
In particular, $b_1, \ldots, b_m, b_{m+k+1}, \ldots, b_n \notin \bar{a}$.
Since $\bar{a}, \bar{b} \in X$ we have
$\tp_\mcM(\bar{b}) = \tp_\mcM(\bar{a})$ so 
(by $\omega$-categoricity) there is $\bar{c} \in X$ such that
$\tp_\mcM(\bar{b}, \bar{c}) = \tp_\mcM(\bar{a}, \bar{b})$.
Then 
\[
\bar{b} \cap \bar{c} = \{b_1, \ldots, b_k\}   \ 
\text{ and $b_i = c_{m+\gamma(i)}$ for all $i = 1, \ldots ,k$.}
\]
By Observation~\ref{facts about independence in degenerate structures}~(ii), we may assume that 
$c_1, \ldots, c_m, c_{m+k+1}, \ldots, c_n \notin \bar{a} \cup \bar{b}$.
It follows that
\[
\bar{a} \cap \bar{c} \ \subseteq \ \{c_{m+1}, \ldots, c_{m+k}\} \ = \ \{b_1, \ldots, b_k\}.
\]
If $i \leq m$ then (as we concluded above) $b_i \notin \bar{a}$ and therefore we get
\[
\bar{a} \cap \bar{c} \ \subseteq \ \{ b_{m+1}, \ldots, b_k\}.
\]
Since $\bar{a} \cap \bar{b} = \{b_{m+1}, \ldots, b_{m+k}\}$ and $m \geq 1$ it follows that
$\bar{a} \cap \bar{c}$ is a proper subset of $\bar{a} \cap \bar{b}$.
As $\tp_\mcM(\bar{b}, \bar{c}) = \tp_\mcM(\bar{a}, \bar{b})$ we also have  $E(\bar{b}, \bar{c})$.
By transitivity of $E$ we get $E(\bar{a}, \bar{c})$.
\hfill $\square$
\\

\noindent
Let $k$ be minimal such that there are $\bar{a}, \bar{b} \in X$ such that $E(\bar{a}, \bar{b})$ and
$|\bar{a} \cap \bar{b}| = k$. By Lemma~\ref{second lemma on disjoint tuples for omega-categorical}, $k > 0$.

\begin{lem}\label{getting I for tuples with minimal intersection for omega-categorical}
There is $I \subseteq \{1, \ldots, n\}$ such that $|I| = k$ and for all $\bar{a}, \bar{b} \in X$
such that $E(\bar{a}, \bar{b})$ and
$|\bar{a} \cap \bar{b}| = k$, we have $\bar{a} \cap \bar{b} = \{a_i : i \in I\} = \{b_i : i \in I\}$.
\end{lem}

\noindent
{\bf Proof.}
For every $\bar{a} \in X$
let $k_{\bar{a}}$ be minimal such that there is $\bar{b} \in X$ such that
$E(\bar{a}, \bar{b})$ and $|\bar{a} \cap \bar{b}| = k_{\bar{a}}$. 
Since all $\bar{a} \in X$ have the same complete type it follows that 
for all $\bar{a} \in X$ there is $\bar{b} \in X$ such that $E(\bar{a}, \bar{b})$
and $|\bar{a} \cap \bar{b}| = k$.
Hence $k_{\bar{a}} = k$ for all $\bar{a} \in X$.

By the minimality of $k$ and 
Lemma~\ref{the most technical lemma about intersection of tuples for omega-categorical},
if  $\bar{a}, \bar{b} \in X$, $E(\bar{a}, \bar{b})$ and 
$|\bar{a} \cap \bar{b}| = k$, then there is 
$I_{\bar{a}, \bar{b}} \subseteq \{1, \ldots, n\}$ such that $|I_{\bar{a}, \bar{b}}| = k$ and
 $\bar{a} \cap \bar{b} = \{a_i : i \in I_{\bar{a}, \bar{b}}\} = \{b_i : i \in I_{\bar{a}, \bar{b}}\}$.
Suppose that for some other $\bar{c} \in X$ we have $E(\bar{a}, \bar{c})$,
$|\bar{a} \cap \bar{c}| = k$ and $I_{\bar{a}, \bar{c}} \neq I_{\bar{a}, \bar{b}}$.
By Observation~\ref{facts about independence in degenerate structures}~(ii),
we may assume that for every $i \notin I_{\bar{a}, \bar{c}}$,
$c_i \notin \bar{b}$.
It follows that $|\bar{b} \cap \bar{c}| < k$ and, by transitivity of $E$,
that $E(\bar{b}, \bar{c})$, which contradicts the choice of $k$.
Hence we conclude that $I_{\bar{a}, \bar{b}} = I_{\bar{a}, \bar{c}}$ for all $\bar{b}, \bar{c} \in X$.
Thus we denote $I_{\bar{a}, \bar{b}}$ by $I_{\bar{a}}$ for any $\bar{b} \in X$.
As all $\bar{a} \in X$ have the same complete type we have $I_{\bar{a}} = I_{\bar{b}}$ for all $\bar{a}, \bar{b} \in X$.
So we denote $I_{\bar{a}}$ by $I$ for any $\bar{a} \in X$.
\hfill $\square$
\\

\noindent
Let $I$ be as in Lemma~\ref{getting I for tuples with minimal intersection for omega-categorical}.
To simplify notation and without loss of generality we assume that
\[
I = \{1, \ldots, k\}.
\]
For any $\bar{a} \in X$ let $\Gamma$ be the set permutations $\gamma$ of $I$ such that
for some $\bar{b} \in X$, $E(\bar{a}, \bar{b})$, $\bar{a} \cap \bar{b} = \{a_i : i \in I\}$ and $a_i = b_{\gamma(i)}$
for all $i \in I$. As all $\bar{a} \in X$ have the same complete type, $\Gamma$ does not depend on $\bar{a}$.
By the transitivity of $E$ and since all $\bar{a} \in X$ have the same complete type 
it follows that $\Gamma$ is closed under composition.
By the symmetry of $E$, $\Gamma$ is closed under inverses.
Hence $\Gamma$ is a {\em group} of permutations of $I$.

\begin{lem}\label{if E holds then there is a permutation for omega-categorical}
If $\bar{a}, \bar{b} \in X$ and $E(\bar{a}, \bar{b})$
then there is $\gamma \in \Gamma$ such that, for all $i \in I$, $a_i = b_{\gamma(i)}$.
\end{lem}

\noindent
{\bf Proof.}
Suppose that $\bar{a}, \bar{b} \in X$ and $E(\bar{a}, \bar{b})$.
Since all tuples in $X$ have the same complete type and $\Gamma$ contains the identity permutation
there is $\bar{c} \in X$ such that
\begin{align*}
E(\bar{b}, \bar{c}), \
\bar{b} \cap \bar{c} = \{b_i : i \in I\} \ \text{ and $c_i = b_i$ for all $i \in I$}.
\end{align*}
Moreover (by Observation~\ref{facts about independence in degenerate structures}~(ii)), we may assume that 
$\bar{a} \cap \bar{c} \subseteq \{c_i : i \in I\}$.
By transitivity of $E$ we have $E(\bar{a}, \bar{c})$.
Therefore $\bar{a} \cap \bar{c} = \{c_i : i \in I\}$ by choice of $k$.
Suppose that
\[
\{a_i : i \in I\} \neq \bar{a} \cap \bar{c}.
\]
It follows that $\{i : a_i \in \bar{c}\} \neq \{i : c_i \in \bar{a}\}$.
Then Lemma~\ref{the most technical lemma about intersection of tuples for omega-categorical}
implies that there is $\bar{d} \in X$ such that $\bar{c} \cap \bar{d} \subsetneq \bar{c} \cap \bar{a}$ and $E(\bar{c}, \bar{d})$,
which contradicts the choice of $k$.
Thus we conclude that $\{a_i : i \in I\} = \bar{a} \cap \bar{c} = \{c_i : i \in I\}$,
and hence there is 
$\gamma \in \Gamma$ such
that $a_i = c_{\gamma(i)}$ for all $i \in I$.
As $b_i = c_i$ for all $i \in I$ we get $a_i = b_{\gamma(i)}$ for all $i \in I$.
\hfill $\square$
\\

\noindent
From Lemma~\ref{if E holds then there is a permutation for omega-categorical}
it follows that if $E'$ is defined as in~(b) of 
Theorem~\ref{characterization of equivalence relations for omega-categorical},
then $E \subseteq E'$.
For the rest of this section let $E''$ be defined as in~(c) of 
Theorem~\ref{characterization of equivalence relations for omega-categorical}.
It remains to prove that $E'' \subseteq E$.

For $\gamma \in \Gamma$ and $\bar{a} = (a_1,  \ldots, a_k)$ we use the notation
$\gamma(\bar{a}) = (a_{\gamma(1)}, \ldots, a_{\gamma(k)})$.

\begin{lem}\label{getting tuples satisfying E with the same type over acl(bar-a')}
There are $\bar{a}' \in M^k$, $\bar{a}^*, \bar{b}^* \in M^{n-k}$ such that
$\bar{a}'\bar{a}^*, \bar{a}'\bar{b}^* \in X$, $E(\bar{a}'\bar{a}^*, \bar{a}'\bar{b}^*)$,
$\bar{a}^* \cap \bar{b}^* = \es$
and $\tp_{\mcM\meq}(\bar{a}^* / \acl_{\mcM\meq}(\bar{a}')) = \tp_{\mcM\meq}(\bar{b}^* / \acl_{\mcM\meq}(\bar{a}'))$.
\end{lem}

\noindent
{\bf Proof.}
By the choice of $k$ and $I$ 
(and since the identity on $I$ belongs to $\Gamma$)
there are $\bar{a}' \in M^k$, $\bar{b}', \bar{c}' \in M^{n-k}$ such that
$\bar{a}'\bar{b}', \bar{a}'\bar{c}' \in X$, $E(\bar{a}'\bar{b}', \bar{a}'\bar{c}')$ and
$\bar{b}' \cap \bar{c}' = \es$.
By Observation~\ref{facts about independence in degenerate structures}~(ii), there are in fact
$\bar{c}'_i \in M^{n-k}$ for all $i < \omega$ such that 
$\bar{a}'\bar{c}'_i \in X$, $E(\bar{a}'\bar{b}', \bar{a}'\bar{c}'_i)$,
$\bar{b}' \cap \bar{c}'_i = \es$ and
$\bar{c}'_i \cap \bar{c}'_j = \es$ for all $i < j < \omega$.
By $\omega$-categoricity there must be $i < j < \omega$ such that
\[
\tp_{\mcM\meq}(\bar{c}'_i /  \acl_{\mcM\meq}(\bar{a}')) = \tp_{\mcM\meq}(\bar{c}'_j / \acl_{\mcM\meq}(\bar{a}')).
\]
By symmetry and transitivity we get $E(\bar{a}'\bar{c}'_i, \bar{a}'\bar{c}'_j)$, so we
are done by taking $\bar{a}^* = \bar{c}'_i$ and $\bar{b}^* = \bar{c}'_j$.
\hfill $\square$

\begin{lem}\label{characterization of E for tuples with minimal intersection for omega-categorical}
Let $\bar{a}' \in M^k$, $\bar{a}^*, \bar{b}^* \in M^{n-k}$ and suppose that
$\bar{a}'\bar{a}^*, \bar{a}'\bar{b}^* \in X$,
$E''(\bar{a}'\bar{a}^*,  \bar{a}'\bar{b}^*)$ 
and $\bar{a}^* \cap \bar{b}^* = \es$.
Then $E(\bar{a}'\bar{a}^*, \bar{a}'\bar{b}^*)$.
\end{lem}

\noindent
{\bf Proof.}
Let $\bar{a}', \bar{a}^*$ and $\bar{b}^*$ satisfy the assumptions of the lemma.
Then, by the definition of $E''$,
\[
\tp_{\mcM\meq}(\bar{a}^* / \acl_{\mcM\meq}(\bar{a}')) = \tp_{\mcM\meq}(\bar{b}^* / \acl_{\mcM\meq}(\bar{a}')).
\]

Since all tuples in $X$ have the same type over $\es$ (in fact even over $\acl_{\mcM\meq}(\es)$)
it follows from
Lemma~\ref{getting tuples satisfying E with the same type over acl(bar-a')}
that there are $\bar{c}^*, \bar{d}^* \in M^{n-k}$ such that
$\bar{a}'\bar{c}^*, \bar{a}'\bar{d}^* \in X$,
$E(\bar{a}'\bar{a}^*, \bar{a}'\bar{c}^*)$, $E(\bar{a}'\bar{b}^*, \bar{a}'\bar{d}^*)$,
$\bar{c}^* \cap \bar{a}^* = \es$, $\bar{d}^* \cap \bar{b}^* = \es$,
\begin{align*}
\tp_{\mcM\meq}(\bar{c}^* / \acl_{\mcM\meq}(\bar{a}')) &= \tp_{\mcM\meq}(\bar{a}^* / \acl_{\mcM\meq}(\bar{a}')) 
\ \text{ and } \\
\tp_{\mcM\meq}(\bar{d}^* / \acl_{\mcM\meq}(\bar{a}')) &= \tp_{\mcM\meq}(\bar{b}^* / \acl_{\mcM\meq}(\bar{a}')).
\end{align*}
Then
\[
\tp_{\mcM\meq}(\bar{c}^* / \acl_{\mcM\meq}(\bar{a}')) = \tp_{\mcM\meq}(\bar{d}^* / \acl_{\mcM\meq}(\bar{a}')).
\]
By Observation~\ref{facts about independence in degenerate structures}~(i), we get
$\bar{a}^* \underset{\bar{a}'}{\ind} \bar{b}^*$, 
$\bar{c}^* \underset{\bar{a}'}{\ind} \bar{a}^*$ and 
$\bar{d}^* \underset{\bar{a}'}{\ind} \bar{b}^*$.
By the independence theorem, there is $\bar{e} \in M^{n-k}$ such that
\begin{align*}
\tp_{\mcM\meq}(\bar{e}^*, \bar{a}^* / \acl_{\mcM\meq}(\bar{a}')) &= 
\tp_{\mcM\meq}(\bar{c}^*, \bar{a}^* / \acl_{\mcM\meq}(\bar{a}')) 
\ \text{ and } \\
\tp_{\mcM\meq}(\bar{e}^*, \bar{b}^* / \acl_{\mcM\meq}(\bar{a}')) &= 
\tp_{\mcM\meq}(\bar{d}^*, \bar{b}^* / \acl_{\mcM\meq}(\bar{a}')).
\end{align*}
This implies that $\bar{a}'\bar{e}^* \in X$, $E(\bar{a}'\bar{a}^*, \bar{a}'\bar{e}^*)$ and 
$E(\bar{a}'\bar{b}^*, \bar{a}'\bar{e}^*)$, so by symmetry and transitivity of $E$
we get $E(\bar{a}'\bar{a}^*, \bar{a}'\bar{b}^*)$.
\hfill $\square$

\begin{lem}\label{characterization of E in general for omega-categorical}
Let $\bar{a}' \in M^k$, $\bar{a}^*, \bar{b}^* \in M^{n-k}$ and suppose that
$\bar{a}'\bar{a}^*, \bar{a}'\bar{b}^* \in X$ and
$E''(\bar{a}'\bar{a}^*,  \bar{a}'\bar{b}^*)$.
Then $E(\bar{a}'\bar{a}^*, \bar{a}'\bar{b}^*)$.
\end{lem}

\noindent
{\bf Proof.}
From $E''(\bar{a}'\bar{a}^*,  \bar{a}'\bar{b}^*)$ we get
\[
\tp_{\mcM\meq}(\bar{a}^* / \acl_{\mcM\meq}(\bar{a}')) = \tp_{\mcM\meq}(\bar{b}^* / \acl_{\mcM\meq}(\bar{a}')).
\]
By Observation~\ref{facts about independence in degenerate structures}~(ii),
there is $\bar{c}^*$ such that $\bar{a}'\bar{c}^* \in X$ and
\[
\tp_{\mcM\meq}(\bar{c}^* / \acl_{\mcM\meq}(\bar{a}')) = \tp_{\mcM\meq}(\bar{b}^* / \acl_{\mcM\meq}(\bar{a}'))
\]
and $\bar{c}^* \cap (\bar{a}^* \cup \bar{b}^*) = \es$.
From the definition of $E''$ we get $E''(\bar{a}'\bar{a}^*, \bar{a}'\bar{c}^*)$ and
$E''(\bar{a}'\bar{c}^*, \bar{a}'\bar{b}^*)$.
Since $\bar{a}^* \cap \bar{c}^* = \es$ and $\bar{b}^* \cap \bar{c}^* = \es$ we get
$E(\bar{a}'\bar{a}^*, \bar{a}'\bar{c}^*)$ and
$E(\bar{a}'\bar{c}^*, \bar{a}'\bar{b}^*)$ by
Lemma~\ref{characterization of E for tuples with minimal intersection for omega-categorical}.
Hence $E(\bar{a}'\bar{a}^*, \bar{a}'\bar{b}^*)$ by the transitivity of $E$.
\hfill $\square$
\\

\noindent
Now 
Theorem~\ref{characterization of equivalence relations for omega-categorical}
follows from 
Lemmas~\ref{if E holds then there is a permutation for omega-categorical}
and~\ref{characterization of E in general for omega-categorical}.

\section{A variation of Theorem~\ref{characterization of equivalence relations for omega-categorical} 
and an application to the generic tetrahedron-free 3-hypergraph}
\label{A variation on the theme}

\noindent
{\bf \em Throughout this section, we assume that $\mcM$ is $\omega$-categorical and that
$\acl_\mcM$ is degenerate.}
Let $0 < n < \omega$, let $p(\bar{x}) \in S^\mcM_n(\es)$ and let 
\[
X =  \{\bar{a} \in M^n : \mcM \models p(\bar{a}) \}.
\]
Suppose that $E$ is a $\es$-definable nontrivial equivalence relation on $X$.
Also assume the following:
\begin{itemize}
\item[$(*)$] {\em If $\bar{a}', \bar{a}^*, \bar{b}^*, \bar{c}^*, \bar{d}^* \in M$ are tuples
such that their ranges are mutually disjoint, the starred sequences are nonemtpy (but $\bar{a}'$ is allowed to be empty)
and $\bar{a}'\bar{a}^*, \bar{a}'\bar{b}^*, \bar{a}'\bar{c}^*$, \\ 
$\bar{a}'\bar{d}^* \in X$,
then there is $\bar{e}^* \in M$ such that 
\[
\tp_{\mcM}(\bar{e}^*, \bar{a}^*, \bar{a}') = \tp_{\mcM}(\bar{c}^*, \bar{a}^*, \bar{a}')
\ \text{ and } \
\tp_{\mcM}(\bar{e}^*, \bar{b}^*, \bar{a}') = \tp_{\mcM}(\bar{d}^*, \bar{b}^*, \bar{a}').
\]
}\end{itemize}

\noindent
Let $E'$ be defined as in part~(b) of
Theorem~\ref{characterization of equivalence relations for omega-categorical}.
Then Lemma~\ref{first lemma on disjoint tuples for omega-categorical} is proved with the help of~$(*)$
(and empty $\bar{a}'$); there is no need to use the independence theorem or the notions of simplicity or SU-rank.
The other lemmas in
Section~\ref{Definable equivalence relations}
up to (and including)
Lemma~\ref{if E holds then there is a permutation for omega-categorical}
are proved in the same way as before.

Define $E''$ on $X$ by
\[
E''(\bar{a}, \bar{b}) \text{ if and only if $a_i = b_i$ for all $i \in I$}.
\]
Then $E''$ is an equivalence relation which is $\es$-definable in $\mcM$.
For this $E''$, the statement of 
Lemma~\ref{characterization of E for tuples with minimal intersection for omega-categorical}
is proved similarly as before, but using~$(*)$ instead of the independence theorem and all
occurences of  `$\acl_{\mcM\meq}(\bar{a}')$' are replaced by `$\bar{a}'$'.
(The same substitutions can be made in 
Lemma~\ref{getting tuples satisfying E with the same type over acl(bar-a')}
and its proof.)
Lemma~\ref{characterization of E in general for omega-categorical}
follows from Lemma~\ref{characterization of E for tuples with minimal intersection for omega-categorical}
in exactly the same way for $E''$ as defined in this section
as was the case in 
Section~\ref{Definable equivalence relations}.
It follows that $E'' \subseteq E \subseteq E'$, where $E'$ is defined in terms of some group $\Gamma$ of permutations of $I$.
Let $\Sigma$ be the set of all $\gamma \in \Gamma$ such that there are $\bar{a}, \bar{b} \in X$ such that
$E(\bar{a}, \bar{b})$ and $a_i =  b_{\gamma(i)}$ for all $i \in I$.
Since $E'' \subseteq E$ the identity function on $I$ belongs to $\Sigma$.
As $E$ is symmetric and transitive it follows that $\Sigma$ is closed under inverses and compositions.
Hence $\Sigma$ is a group of permutations.
Thus we get the following version of
Theorem~\ref{characterization of equivalence relations for omega-categorical}:

\begin{theor}\label{characterization of equivalence relations in a special case}
Suppose that $\mcM$ is $\omega$-categorical and that $\acl_\mcM$ is degenerate.
Moreover suppose that~$(*)$ holds for $X$ as defined in this section and that $E$ is a nontrivial $\es$-definable 
equivalence relation on $X$.
Then there is a nonempty $I \subseteq \{1, \ldots, n\}$ and a group $\Gamma$ of permutations of $I$ such that 
for all $\bar{a} = (a_1, \ldots, a_n), \bar{b} = (b_1, \ldots, b_n) \in X$, $E(\bar{a}, \bar{b})$ if and only if
there is $\gamma \in \Gamma$ such that $a_i = b_{\gamma(i)}$ for all $i \in I$.
\end{theor}

\begin{cor}\label{algebraic closure equals definable closure}
Suppose that $\mcM$ is $\omega$-categorical  and that $\acl_\mcM$ is degenerate.
Moreover, assume that~$(*)$ holds for any choice of $0 < n < \omega$ and any $p \in S^\mcM_n(\es)$.
Then, for every $A \subseteq M$, $\acl_{\mcM\meq}(A) = \dcl_{\mcM\meq}(A)$.
\end{cor}

\noindent
{\em Remark:} For the rest of this section we use the following notation.
For every $n < \omega$, every $\es$-definable equivalence relation $E$ on $M^n$ and every 
$\bar{a} \in M^n$, let $[\bar{a}]_E$ denote the $E$-equivalence class of $\bar{a}$ {\em as an element of $M\meq$}
(and not as a subset of $M^n$).

\medskip

\noindent
{\bf Proof.}
Suppose that $A \subseteq M$, $b \in M\meq$ and $b \in \acl_{\mcM\meq}(A)$.
Without loss of generality we may assume that $A$ is finite.
For some $n < \omega$, some $\es$-definable equivalence relation $E$ on $M^n$ and
some $\bar{b} = (b_1, \ldots, b_n) \in M^n$ we have 
$b = [\bar{b}]_E$.

Let $p(\bar{x}) = \tp(\bar{b})$ and $X = \{\bar{a} \in M^n : \mcM \models p(\bar{a})\}$.
Without loss of generality we may assume that $p(\bar{x})$ implies $x_i \neq x_j$ for all $1 \leq i < j \leq n$.\footnote{
Because otherwise $b$ is interdefinable with some $b'$ for which this assumption holds.}
If for all $\bar{a}, \bar{a}' \in X$, $E(\bar{a}, \bar{a}')$, then 
$[\bar{b}]_E \in \dcl_{\mcM\meq}(\es)$.
If for all $\bar{a} \in X$ such that $\bar{a} \neq \bar{b}$ we have $\neg E(\bar{a}, \bar{b})$,
then $\bar{b} \in \acl_{\mcM\meq}([\bar{b}]_E) \subseteq \acl_{\mcM\meq}(A)$ and hence
$\bar{b} \in \acl_\mcM(A) = A$. It follows that $[\bar{b}]_E \in \dcl_{\mcM\meq}(A)$.

Hence, from now on we can assume that $E$ is nontrivial on $X$.
By Theorem~\ref{characterization of equivalence relations in a special case},
there is a nonempty $I \subseteq \{1, \ldots, n\}$ and a group $\Gamma$ of permutations of $I$ such that
for all $\bar{a} = (a_1, \ldots, a_n),  \bar{a}' = (a'_1, \ldots, a'_n) \in X$, $E(\bar{a}', \bar{a}')$ 
and only if for some $\gamma \in \Gamma$ and all $i \in I$, $a_i = a'_{\gamma(i)}$.
Without loss of generality assume that $I = \{1, \ldots, k\}$ where $k \leq n$.

First suppose that $b_1, \ldots, b_k \in A$. 
Let $\psi(x, b_1, \ldots, b_k)$ be the formula (with no other parameters than $b_1, \ldots, b_k$)
which expresses:
\begin{itemize}
\item[] ``$x$ has sort $E$ and if $x = [(y_1, \ldots, y_n)]_E$ then there is a permutation $\gamma \in \Gamma$
such that for all $i = 1, \ldots, k$, $b_i = y_{\gamma(i)}$.''
\end{itemize}
Clearly, $b = [\bar{b}]_E$ satisfies $\psi(x, b_1, \ldots, b_k)$.
And if $b' = [(b'_1, \ldots, b'_n)]_E$ satisfies $\psi(x, b_1, \ldots, b_k)$,
then for some permutation $\gamma \in \Gamma$, $b_i = b'_{\gamma(i)}$ for all $i = 1, \ldots, k$
and hence $E(\bar{b}, (b'_1, \ldots, b'_n))$, so $b = b'$.
Thus $b \in \dcl_{\mcM\meq}(b_1, \ldots, b_k) \subseteq \dcl_{\mcM\meq}(A)$.

Now suppose that for some $i \in I = \{1, \ldots, k\}$, $b_i \notin A$ (and we will derive a contradiction from this).
To simplify notation, and without loss of generality, assume that $i = 1$.
Then $b_1 \notin \acl_{\mcM\meq}(A \cup \{b_2, \ldots, b_n\})$ 
so for every $j < \omega$ there is
$b_{1, j} \in M \setminus (A \cup \{b_{1, l} : l < j\})$ such that
\[
\tp_\mcM(b_{1, j}, b_2, \ldots, b_n / A) = \tp_\mcM(b_1, b_2, \ldots, b_n / A).
\]
Let
$\bar{b}^j = (b_{1, j}, b_2, \ldots, b_n)$ for all $j < \omega$.
Then for every  $j < \omega$ there is an automorphism of $\mcM\meq$ which sends $\bar{b}$ to $\bar{b}^j$ and fixes $A$ pointwise.
As this automorphism sends $[\bar{b}]_E$ to $[\bar{b}^j]_E$ we get
$\tp_{\mcM\meq}([\bar{b}^j]_E / A) = \tp_{\mcM\meq}([\bar{b}]_E / A)$ for all $j < \omega$.
By the choice of $b_{1, j}$ we have $[\bar{b}^j]_E \neq [\bar{b}^l]_E$ whenever $j \neq l$,
so this contradicts the assumption that $b = [\bar{b}]_E \in \acl_{\mcM\meq}(A)$.
\hfill $\square$

\begin{rem}\label{the results about equivalence relations hold for all models of the theory}{\rm
Note that Theorems~\ref{characterization of equivalence relations for omega-categorical},
\ref{characterization of equivalence relations in a special case}
and Corollary~\ref{algebraic closure equals definable closure}
are really results about $Th(\mcM)$ (rather than just $\mcM$). 
More precisely, the characterizations of
equivalence relations in 
Theorems~\ref{characterization of equivalence relations for omega-categorical}
and~\ref{characterization of equivalence relations in a special case} 
hold if $\mcM$ (in the respective theorem) is replaced by any model of $Th(\mcM)$.
It follows (from its proof) that $\mcM$ in Corollary~\ref{algebraic closure equals definable closure}
can be replaced by any model of $Th(\mcM)$.
}\end{rem}

\begin{exam}\label{examples satisfying the special case}{\rm
(i) It is easy to see that if $\mcM$ is a binary random structure, then~$(*)$ is satisfied for
any choice of $0 < n < \omega$ and any $p \in S^\mcM_n(\es)$. \\
(ii) Let $\mcM$ be the generic tetrahedron-free 3-hypergraph
given by Definition~\ref{definition of generic tetrahedron-free 3-hypergraph}. 
It is straightforward to verify that $\mcM$ satisfies~$(*)$ for any choice of $0 < n < \omega$ and
any $p \in S^\mcM_n(\es)$. The reason is essentially that (with the notation of~$(*)$) we can find $\bar{e}^*$
such that
\[
\tp_{\mcM}(\bar{e}^*, \bar{a}^*, \bar{a}') = \tp_{\mcM}(\bar{c}^*, \bar{a}^*, \bar{a}')
\ \text{ and } \
\tp_{\mcM}(\bar{e}^*, \bar{b}^*, \bar{a}') = \tp_{\mcM}(\bar{d}^*, \bar{b}^*, \bar{a}').
\]
and if $e \in \rng(\bar{e}^*)$, $a \in \rng(\bar{a}^*)$ and $b \in \rng(\bar{b}^*)$ then
$\{e, a, b\}$ is {\em not} a hyperedge.
}\end{exam}

\noindent
Recall that Remark~\ref{remark about the generic tetrahedron-free 3-hypergraph}
lists some known properties of the generic tetrahedron-free 3-hypergraph.
The following result implies that it is also 1-based\footnote{The fact that the generic tetrahedron-free 3-hypergraph is 1-based
also follows from a result by Conant in \cite{Con}.}:

\begin{prop}\label{the generic tetrahedron-free hypergraph is 1-based}
Suppose that $\mcM$ is countable, homogeneous and supersimple with SU-rank 1.
Moreover, assume that $\acl_\mcM$ is degenerate and that~$(*)$ holds.
Then $\mcM$ is 1-based.
\end{prop}

\noindent
{\bf Proof.}
For any type $p$ and structure $\mcN$, let $\real_\mcN(p)$ be the set of tuples of elements from $\mcN$ which realize $p$.
By \cite[Corollary~4.7]{HKP}, where 1-based theories are called ``modular'',
it suffices to prove the following:
\begin{itemize}
\item[($\dagger$)] If $\mcN \models Th(\mcM)$, $A \subseteq N\meq$ and $p(x)$ is a complete type over $A$ (possibly realized
by imaginary elements) with SU-rank 1,
then  $(\real_{\mcN\meq}(p), \cl)$, where $\cl(B) = \acl_{\mcN\meq}(B \cup A) \cap \real_{\mcN\meq}(p)$
for all $B \subseteq \real_{\mcN\meq}(p)$, is a trivial pregeometry
(i.e. if $a \in \cl(B)$ then $a \in \cl(b)$ for some $b \in B$).\footnote{
If $(\real_{\mcN\meq}(p), \cl)$ is a trivial pregeometry then $p$ is ``modular'' in the sense of \cite{HKP}.}
\end{itemize}
The first step is to show that it suffices to prove~($\dagger$) in the case when $A$ is finite.

Suppose that $\mcN \models Th(\mcM)$, $A \subseteq N\meq$ and that $p(x)$ is a complete type over $A$ with SU-rank 1.
Moreover, suppose that $(\real_{\mcN\meq}(p), \cl)$ is nontrivial. So there are finite
$B \subseteq \real_{\mcN\meq}(p)$ and $a \in \real_{\mcN\meq}(p)$ such that
$a \in \acl_{\mcN\meq}(B  \cup A) \setminus \acl_{\mcN\meq}(\{b\} \cup A)$ for all $b \in B$.
Then there is finite $A_0 \subseteq A$ such that 
$a \in \acl_{\mcN\meq}(B  \cup A_0) \setminus \acl_{\mcN\meq}(\{b\} \cup A_0)$ for all $b \in B$.
Since $Th(\mcM)$ is supersimple there is a finite $A_1 \subseteq A$ such that $p$ does not divide over $A_1$.
Let $p'$ be the restriction of $p$ to formulas with parameters from $A' = A_0 \cup A_1$. 
Then $p'$ has SU-rank 1 and $a \in \acl_{\mcN\meq}(B  \cup A') \setminus \acl_{\mcN\meq}(\{b\} \cup A')$ for all $b \in B$,
so $(\real_{\mcN\meq}(p'), \cl')$ is a nontrivial pregeometry,
where $\cl'(B) = \acl_{\mcN\meq}(B \cup A') \cap \real_{\mcN\meq}(p')$
for all $B \subseteq \real_{\mcN\meq}(p')$.
Hence it suffices to prove~($\dagger$) for finite $A$.

The next step is to show that it suffices to prove~($\dagger$) in the case when $A$ is finite {\em and}
a subset of $N$ (so that only ``real elements'' occur in $A$).
Let $A \subseteq N\meq$  be finite and suppose that $(\real_{\mcN\meq}(p), \cl)$ is nontrivial 
(where $p(x)$ is a complete type over $A$ with SU-rank 1).
Suppose that $a, \bar{b} \in \real_{\mcN\meq}(p)$ are such that 
$a \in \acl_{\mcN\meq}(\rng(\bar{b})  \cup A) \setminus \acl_{\mcN\meq}(\{b\} \cup A)$ for all $b \in \rng(\bar{b})$.
There is finite  $C \subseteq N$ such that $A \subseteq \dcl_{\mcN\meq}(C)$.
Let $\bar{c}$ enumerate $C$. By considering a realization of a nondividing extension of 
$\tp_{\mcN\meq}(\bar{c} /A)$ to $A \cup \{a\} \cup \rng(\bar{b})$ we may assume that $a\bar{b} \underset{A}{\ind} C$.
Since $A \subseteq \dcl_{\mcN\meq}(C)$ we have $a \in \acl_{\mcN\meq}(\rng(\bar{b}) \cup C)$.
Let $b \in \rng(\bar{b})$. 
As $\su(p) = 1$ and $a \notin \acl_{\mcN\meq}(\{b\} \cup A)$ we have $a \underset{A}{\ind} b$ and hence
$\su(a,b /A) = 2$. By the choice of $C$ we have $ab \underset{A}{\ind} C$ and therefore $\su(a,b / C) = 2$
which implies that $a \underset{C}{\ind} b$.
Hence $a \notin \acl_{\mcM\meq}(\{b\} \cup C)$.

Thus it suffices to prove~($\dagger$) for finite $A \subseteq N$.
In fact, since $\mcM$ is $\omega$-saturated it suffices to prove that if 
$A \subseteq M$ is finite and $p(x)$ is a complete type over $A$ with SU-rank 1,
then  $(\real_{\mcM\meq}(p), \cl)$ is a trivial pregeometry.
This is the last step of the proof.

Let $A \subseteq M$ be finite and let $p(x)$ be a complete type over $A$ with SU-rank 1.
Suppose, towards a contradiction, that $(\real_{\mcM\meq}(p), \cl)$
is a nontrivial pregeometry.
Then there are $1 < m < \omega$ and distinct $a, b_1, \ldots, b_m \in \real_{\mcM\meq}(p)$ such that
\begin{align}
&a \in \acl_{\mcM\meq}(\{b_1, \ldots, b_m\} \cup A) \ \text{ and} \label{a is in acl of all b-i's} \\
&a \notin \acl_{\mcM\meq}(\{b_k\} \cup A) \ \text{ for every } k = 1, \ldots, m. \label{a is not in acl of b-i}
\end{align}
Let $0 < n < \omega$ and $E^+ \subseteq M^{2n}$ be a $\es$-definable equivalence relation such that
all elements that realize $p$ are of sort $E^+$.
Furthermore, there is
$q \in S^\mcM_n(\es)$ such that for every $a \in \real_{\mcM\meq}(p)$ there is $\bar{a} \in M^n$ such that
$a = [\bar{a}]_{E^+}$ and $\mcM \models q(\bar{a})$.
Let 
\[
X = \{\bar{a} : \mcM \models q(\bar{a}) \}
\]
and let $E$ be the restriction of $E^+$ to $X$.
Since $a, b_1, \ldots, b_m$ are distinct elements it follows that $E$ has at least two classes.
So either $E$ is nontrivial or $E$ is the identity relation on $X$. 
In either case, and using Theorem~\ref{characterization of equivalence relations in a special case} in the first case,
there are a nonempty $I \subseteq \{1, \ldots, n\}$ and a group $\Gamma$ of permutations of $I$ 
such that for all $\bar{a} = (a_1, \ldots, a_n), \bar{b} = (b_1, \ldots, b_n) \in X$, $E(\bar{a}, \bar{b})$ if and only if
there is $\gamma \in \Gamma$ such that $a_i = b_{\gamma(i)}$ for all $i \in I$.
For every $k = 1, \ldots, m$, choose any $\bar{b}_k = (b_{k,1}, \ldots, b_{k,n}) \in X$ such that $[\bar{b}_k]_{E^+} = b_k$.
By the characterization of $E$ and Observation~\ref{facts about independence in degenerate structures}~(ii)
it follows that
there is $\bar{a} = (a_1, \ldots, a_n) \in X$ such that 
\begin{align*}
&a = [\bar{a}]_{E^+}, \\
&\{a_i : i \notin I\} \cap \{b_{k,i} : i \notin I \} = \es \ \text{ for all $k = 1, \ldots, m$ and } \\
&\{a_i : i \notin I\} \cap A = \es.
\end{align*}
We now divide the argument into two cases, both of which will lead to contradictions.

\medskip

\noindent
{\bf \em Case 1.} Suppose that for every $k \in \{1, \ldots, m\}$,
\[
\big( \{a_i : i \in I\} \setminus A \big) \cap \big( \{b_{k,i} : i \in I \} \setminus A \big) = \es.
\]

\noindent
By Observation~\ref{facts about independence in degenerate structures}~(i),
$\bar{a} \underset{A}{\ind}\{b_{k,i} : k = 1, \ldots, m, \ i = 1, \ldots, n \}$
which (since $a \in \acl_{\mcM\meq}(\bar{a})$ and $b_k \in \acl_{\mcM\meq}(\bar{b}_k)$)
implies that $a \underset{A}{\ind} \{b_1, \ldots, b_m\}$ and hence (since $\mcM$ has SU-rank 1)
\[
a \notin \acl_{\mcM\meq}(\{b_1, \ldots, b_m\} \cup A).
\] This contradicts~(\ref{a is in acl of all b-i's}).

\medskip

\noindent
{\bf \em Case 2.} Suppose that for some $k \in \{1, \ldots, m\}$,
\[
\big( \{a_i : i \in I\} \setminus A \big) \cap \big( \{b_{k,i} : i \in I \} \setminus A \big) \neq \es.
\]

\noindent
Then, by the characterization of $E$, there is $i \in I$ such that 
\begin{itemize}
\item $a_i \notin A$ and
\item if $\bar{a}' \in X$ is such that $[\bar{a}']_{E^+} = a$ and $\bar{b}'_k \in X$ is such that $[\bar{b}'_k]_{E^+} = b_k$,
then $a_i \in \rng(\bar{a}') \cap \rng(\bar{b}'_k)$.
\end{itemize}
Now it is straightforward to prove (for example by using the definition of dividing; details are left to the reader) that
$a \underset{A}{\nind}b_k$ and hence $a \in \acl_{\mcM\meq}(\{b_k\} \cup A)$.
But this contradicts~(\ref{a is not in acl of b-i}).
\hfill $\square$

\begin{rem}\label{proof of an earlier fact}{\rm
Finally we explain why Fact~\ref{binary homogeneous simple structures are 1-based} holds.
Suppose that $\mcM$ is countable, binary, homogeneous and simple,
so $\mcM$ is supersimple with finite SU-rank (by \cite{Kop15}).
We will show that $\mcM$ has trivial dependence.

Consider the following statement for any simple $\omega$-saturated $\mcN$:

\begin{itemize}
\item[($\bigstar$)] if $A \subseteq N$ is finite, $0 < n < \omega$ and $a_1, \ldots, a_n \in N\meq$ are pairwise independent over $A$,
then $\{a_1, \ldots, a_n\}$ is an independent set over $A$.
\end{itemize}
Lemma~1 in \cite{Goode} says that if $\mcN$ is stable and~($\bigstar$) holds in the special case when $a_1, \ldots, a_n \in N$, 
then it also holds in the generality stated above. Its proof uses only basic properties of forking/dividing
which also hold in simple theories/structures, as observed by Palac\'{i}n \cite{Pal}. Therefore we conclude that if $\mcN$ is simple 
and~($\bigstar$) holds in the special case when $a_1, \ldots, a_n \in N$, then it also holds in the generality stated above.

Now suppose that $\mcN \models Th(\mcM)$, so $\mcN$ is $\omega$-saturated (since $Th(\mcM)$ is $\omega$-categorical).
By \cite[Corollary 6]{Kop15}, ($\bigstar$) holds for any $0 < n < \omega$ and any $a_1, \ldots, a_n \in N$.
Hence, ($\bigstar$) holds in the generality stated above.

In order to prove that $\mcM$ has trivial dependence it suffices, 
according to the argument in the beginning of the proof of Proposition~\ref{the generic tetrahedron-free hypergraph is 1-based},
to prove the following:

\begin{itemize}
\item[($\spadesuit$)] Suppose that $\mcN \models Th(\mcM)$, $A \subseteq N$ is finite and $p(x)$ is a complete type 
(possibly realized by imaginary elements) over $A$ with SU-rank 1.
Then  $(\real_{\mcN\meq}(p), \cl)$, where $\cl(B) = \acl_{\mcN\meq}(B \cup A) \cap \real_{\mcN\meq}(p)$
for all $B \subseteq \real_{\mcN\meq}(p)$, is a trivial pregeometry
\end{itemize}

\noindent
But ($\spadesuit$) is a direct consequence of ($\bigstar$), so we are done.\footnote{
It is tempting to claim that one could use Lemmas~1, 4 and Proposition~5 of \cite{Goode} to
directly conclude that $\mcM$ has trivial dependence. But the proof of Proposition~5 refers to the proof
Proposition~2 (in \cite{Goode}) which uses the notion of ``heir'', which is one reason why it is
not clear to me how the argument would be translated to the context of simple structures.}
}\end{rem}

\noindent
{\bf Acknowledgement.} I would like to thank the anonymous referee for
numerous suggestions of how to improve and clarify the arguments.

\end{document}